\documentclass[11pt, reqno]{amsart}

 \usepackage{amsmath}
 \usepackage{amssymb}
\usepackage{bm}

\DeclareMathAccent{\mathring}{\mathalpha}{operators}{"17}

 \usepackage{color}

\newcommand{\mysection}[1]{\section{#1}
      \setcounter{equation}{0}}

\newcommand{\nlimsup}{\operatornamewithlimits{\overline{lim}}}
\newcommand{\nliminf}{\operatornamewithlimits{\underline{lim}}}

\newtheorem{theorem}{Theorem}[section]
\newtheorem{lemma}[theorem]{Lemma}

\theoremstyle{definition}
\newtheorem{assumption}{Assumption}[section]
\newtheorem{definition}{Definition}[section]
\theoremstyle{remark}
\newtheorem{remark}{Remark}[section]

\newtheorem{example}{Example}[section]

\newcommand{\tr}{\text{\rm tr}\,}

 \makeatletter 
 \def\dashint{%
 \operatorname%
 {\,\,\text{\bf--}\kern-.98em\DOTSI\intop\ilimits@\!\!}}
\def\ninf{\qopname\relax\@empty{inf\phantom{p}\!\!\!}}
 \makeatother

\newcommand\bbeta{\text{\raise-.2ex\hbox{$\bm{\beta}$}}}

\newcommand\bR{\mathbb{R}}

\newcommand\bB{\mathbb{B}}

\newcommand\cF{\mathcal{F}}

\newcommand\cO{\mathcal{O}}

\newcommand\cS{\mathcal{S}}

\newcommand\frA{\mathfrak{A}}
\newcommand\frB{\mathfrak{B}}

\newcommand\frP{\mathfrak{P}}

\newcommand\infsup{\operatornamewithlimits{inf\,\,\,sup}}
\newcommand\supinf{\operatornamewithlimits{sup\,\,\,inf}}

\begin{document}

\title[Dynamic programming principle]
{On the dynamic programming principle for
uniformly nondegenerate stochastic  differential 
games in domains}

\author{N.V. Krylov}
\thanks{The  author was partially supported by
 NSF Grant DNS-1160569}
\email{krylov@math.umn.edu}
\address{127 Vincent Hall, University of Minnesota,
 Minneapolis, MN, 55455}



\begin{abstract}
We prove the dynamic programming principle
for uniformly nondegenerate stochastic differential games
in the framework of time-homogeneous
 diffusion processes considered up to the first exit
time from a domain. The zeroth-order ``coefficient'' and the 
``free'' term 
are only assumed to be measurable.
In contrast with previous results established
for constant stopping times we allow arbitrary
stopping times and randomized ones as well.
The main assumption, which will be removed
in a subsequent article, is that there exists a sufficiently
regular solution of the Isaacs equation.
\end{abstract}

\maketitle

\mysection{Introduction}

The dynamic programming principle is one of the basic tools
in the theory of controlled diffusion processes. In early
70's it allowed
one to obtain results about the unique solvability in 
classes of differentiable
functions of
Bellman's equations, which,  for about ten years, were the only known
results for more or less general fully nonlinear second-order
elliptic and
parabolic equations.

In this paper
we will be only dealing with the dynamic programming principle
for stochastic  differential  games.
 Concerning all other aspect of the theory
of stochastic  differential  games we refer the reader to \cite{BL08},
\cite{FS89}, \cite{Ko09}, \cite{Ni88},
 and \cite{Sw96} and the references therein.

It seems to the author that 
Fleming and  Souganidis in \cite{FS89} were the first authors
who proved the dynamic programming principle  with nonrandom stopping times
for stochastic   differential games in the whole space
on a finite time horizon.
 They used rather involved constructions to overcome
some measure-theoretic
difficulties, a technique somewhat resembling
the one in Nisio \cite{Ni88}, and the theory of viscosity
solutions.

In \cite{Ko09}  Kovats considers time-homogeneous 
stochastic differential games
in a ``weak'' formulation
in smooth {\em domains\/} and proves the dynamic programming principle
again with nonrandom stopping times. He uses approximations
of policies by piece-wise constant ones and proceeds
similarly to \cite{Ni88}.
 
\'Swi{\c e}ch in \cite{Sw96} reverses the arguments
in \cite{FS89} and proves the dynamic programming principle
for time-homogeneous stochastic differential games in the whole space
with constant stopping times ``directly''
from knowing   that the viscosity solutions
exist. His method is quite similar to the
so-called verification principle in the theory
of controlled diffusion processes.

It is also worth mentioning the paper \cite{BL08} by Buckdahn and Li
where the dynamic programming principle
for constant stopping times in the time-inhomogeneous setting
in the whole space is derived by using
the theory of backward forward stochastic equations.

Basically, we adopt the strategy of    \'Swi{\c e}ch
(\cite{Sw96}) which is based on using the fact that in many cases
the Isaacs equation has a sufficiently regular solution.
In \cite{Sw96} viscosity solutions are used and we rely on
classical ones.

The main emphasis of \cite{FS89}, \cite{Ko09}, \cite{Ni88},
 and \cite{Sw96} is on proving that (upper and lower)
  value functions for stochastic  differential  games are viscosity
solutions of the corresponding Isaacs equations
and the dynamic programming principle is
used just as a tool to do that. In our setting 
the zeroth-order coefficient and the running payoff
function  
can be just measurable and in this situation
neither our methods nor 
 the methods based on the notion of viscosity 
solution seem  to be of much help while proving
that the value function is a viscosity solution.

Our main future goal is to develop some tools which would allow us
in a subsequent article to show that the value
functions are of class $C^{0,1}$, provided that the data
are there, for {\em possibly degenerate\/} stochastic
differential games without assuming that the zeroth-order
coefficient is large  enough negative.
 On the way to achieve this goal
one of the main steps, apart from proving the dynamic
programming principle, is to prove certain representation
formulas which will be derived in a subsequent article
from our Theorems \ref{theorem 2.27.2} and \ref{theorem 2.27.3},
 in the first of which the process
is not assumed to be uniformly nondegenerate.
Another important ingredient consists of approximations
results allowing us to approximate stochastic  differential  games
with the ones for which the corresponding
Isaacs equations have sufficiently regular solutions.
This issue will be addressed in a subsequent article.

One of the main results of the present article,
Theorem \ref{theorem 1.14.1}, is about the dynamic
programming principle in a very general form
including stopping and randomized stopping times.
It is proved under the assumption that 
the corresponding
Isaacs equations have sufficiently regular solutions.

In Theorem \ref{theorem 2.18.1} we prove the H\"older
continuity of the value function in our case
where the zeroth-order coefficient and the 
running payoff
function can be discontinuous.

Theorem \ref{theorem 1.14.1} concerns
time-homogeneous stochastic differential games unlike the time inhomogeneous
 in \cite{FS89}
and generalizes the {\em
corresponding\/} results of \cite{Sw96} and \cite{Ko09},
where however  degenerate case is not excluded.

The article is organized as follows.
In Section \ref{section 2.26.3} we state our main results
to which actually, as we have pointed out  implicitly above,
belongs Theorems \ref{theorem 2.27.2} and \ref{theorem 2.27.3} reminding 
the verification principle from the theory of
controlled diffusion processes. The main technical tool 
for proving these theorems is laid out 
 in a rather long Section \ref{section 2.11.1}
for processes which may be degenerate. We prove there
Theorems 
\ref{theorem 2.11.1},
 \ref{theorem 4.17.1}, and  \ref{theorem 4.17.2}.
In a short Section \ref{section 2.27.3} we give their versions
 for uniformly 
nondegenerate case. These versions look stronger but 
Theorem  \ref{theorem 2.27.3}
is proved only for uniformly 
nondegenerate case. In Section \ref{section 3.21.1}
we prove an auxiliary result which allows us to
investigate the behavior of the value function
near the boundary.
In the final short Section \ref{section 4.25.1} we combine
previous results and prove Theorem   \ref{theorem 1.14.1}.

\mysection{Main results}
                                           \label{section 2.26.3}

Let $\bR^{d}=\{x=(x_{1},...,x_{d})\}$
be a $d$-dimensional Euclidean space and  let $d_{1}\geq d$
and $k\geq1$ be  integers.
Assume that we are given separable metric spaces
  $A$, $B$, and $P$  and let,
for each $\alpha\in A$, $\beta\in B$, and $p\in P$
  the following 
  functions on $\bR^{d}$ be given: 

(i) $d\times d_{1}$
matrix-valued $\sigma^{\alpha\beta}(p,x)
=\sigma(\alpha,\beta,p,x)=
(\sigma^{\alpha\beta}_{ij}(p,x))$,

(ii)
$\bR^{d}$-valued $b^{\alpha\beta}(p,x)=
b(\alpha,\beta,p,x)=
(b^{\alpha\beta}_{i }(p,x))$, and

(iii)
real-valued  functions $c^{\alpha\beta}(p,x)=c(\alpha,\beta,p,x)$,   
  $f^{\alpha\beta}(p,x)=f(\alpha,\beta,p,x)$, and  
$g(x)$.

Introduce
$$
a^{\alpha\beta}(p,x):=(1/2)\sigma^{\alpha\beta}(p,x)
(\sigma^{\alpha\beta}(p,x))^{*},
$$
fix a $\bar{p}\in P$, and set
$$
(\bar{\sigma},\bar{a},\bar{b},\bar{c},\bar{f})^{\alpha\beta}(x)
=(\sigma,a,b,c,f) ^{\alpha\beta}( \bar{p},x), 
$$

\begin{assumption}
                                    \label{assumption 1.9.1}
(i) All the above functions are continuous with respect to
$\beta\in B$ for each $(\alpha,p,x)$ and continuous with respect
to $\alpha\in A$ uniformly with respect to $\beta\in B$
for each $(p,x)$. Furthermore, they are Borel measurable in $(p,x)$ for each
$(\alpha,\beta)$, the function
$g(x)$ is bounded and uniformly continuous on $\bR^{d}$, and
$c^{\alpha\beta}\geq0$.

(ii) The functions $\bar{\sigma}^{\alpha\beta}( x )$ and
$\bar{b}^{\alpha\beta}( x )$ are uniformly continuous with
respect to $x$ uniformly with respect to 
$(\alpha,\beta )\in A\times B $
and
for any $x \in\bR^{d}$ 
 and $(\alpha,\beta,p)\in A\times B\times P$
$$
  \|\sigma^{\alpha\beta}(p,x )\|,|b^{\alpha\beta}(p,x )|
\leq K_{0},
$$
where   $K_{0}$ is
a fixed constants and for a matrix $\sigma$ we denote $\|\sigma\|^{2}
=\tr\sigma\sigma^{*}$.
\end{assumption}

Let $(\Omega,\cF,P)$ be a complete probability space,
let $\{\cF_{t},t\geq0\}$ be an increasing filtration  
of $\sigma$-fields $\cF_{t}\subset \cF_{t}$ such that
each $\cF_{t}$ is complete with respect to $\cF,P$, and let
$w_{t},t\geq0$, be a standard $d_{1}$-dimensional Wiener process
given on $\Omega$ such that $w_{t}$ is a Wiener process
relative to the filtration $\{\cF_{t},t\geq0\}$.

 The set of progressively measurable $A$-valued
processes $\alpha_{t}=\alpha_{t}(\omega)$ is denoted by $\frA$. 
Similarly we define $\frB$
as the set of $B$-valued  progressively measurable functions.
By  $ \bB $ we denote
the set of $\frB$-valued functions 
$ \bbeta(\alpha_{\cdot})$ on $\frA$
such that, for any $T\in(0,\infty)$ and any $\alpha^{1}_{\cdot},
\alpha^{2}_{\cdot}\in\frA$ satisfying
\begin{equation}
                                                  \label{4.5.4}
P(  \alpha^{1}_{t}=\alpha^{2}_{t} 
 \quad\text{for almost all}\quad t\leq T)=1,
\end{equation}
we have
$$
P(  \bbeta_{t}(\alpha^{1}_{\cdot})=\bbeta_{t}(\alpha^{2}_{\cdot}) 
\quad\text{for almost all}\quad t\leq T)=1.
$$

 We closely follow the standard setup (but not the notation)
from \cite{FS89}, \cite{Ko09},  
 and \cite{Sw96} allowing the control processes to depend
on the past information contained in $\{\cF_{t}\}$.
By the way, in the situation
of controlled diffusion processes (not stochastic  differential 
 games)
 these control processes were first introduced 
in \cite{Kr72} and turned out to be  extremely useful
in developing  the theory
of Bellman's equations.

\begin{definition}
A 
function $p^{\alpha_{\cdot}\beta_{\cdot}}_{t}=
p^{\alpha_{\cdot}\beta_{\cdot}}_{t}(\omega)$ given
on $\frA\times\frB\times\Omega\times[0,\infty)$
is called {\em a control adapted process} if,
for any $(\alpha_{\cdot},\beta_{\cdot})\in\frA\times\frB$,
it is progressively measurable in $(\omega,t)$ and, for any
$T\in(0,\infty)$, we have
$$
P( p^{\alpha^{1}_{\cdot}\beta^{1}_{\cdot}}_{t}
=p^{\alpha^{2}_{\cdot}\beta^{2}_{\cdot}}_{t}
\quad\text{for almost all}\quad t\leq T)=1
$$
as long as
$$
P(  \alpha^{1}_{t}=\alpha^{2}_{t} ,
\beta^{1}_{t}=\beta^{2}_{t}
\quad\text{for almost all}\quad t\leq T)=1.
$$
The set of control adapted $P$-valued processes is denoted by $\frP$. 
\end{definition}

We fix a $p\in\frP$ (for the rest of the article)
 and for $\alpha_{\cdot}\in\frA$, 
$ \beta_{\cdot} 
\in\frB$, and $x\in\bR^{d}$ consider the following It\^o
equation
\begin{equation}
                                             \label{5.11.1}
x_{t}=x+\int_{0}^{t}\sigma^{\alpha_{s}
\beta_{s} }( p^{\alpha_{\cdot}\beta_{\cdot}}_{s},x_{s})\,dw_{s}
+\int_{0}^{t}b^{\alpha_{s}
\beta_{s} }(  p^{\alpha_{\cdot}\beta_{\cdot}}_{s},x_{s})\,ds.
\end{equation}

\begin{assumption}
                                           \label{assumption 5.23.2}
Equation \eqref{5.11.1} satisfies the usual hypothesis, that is
for any $\alpha_{\cdot}\in\frA$, 
$ \beta_{\cdot} 
\in\frB$, and $x\in\bR^{d}$ it has a unique solution
denoted by $x^{\alpha_{\cdot} 
\beta_{\cdot} x}_{t} $ and $x^{\alpha_{\cdot} 
\beta_{\cdot} x}_{t} $ is a control adapted process
for each $x$.

\end{assumption}

\begin{remark}
                                                      \label{remark 5.23.1}
As is well known, equation \eqref{5.11.1}  satisfies 
the usual hypothesis
if Assumption \ref{assumption 1.9.1} is satisfied and
for any $x,y\in\bR^{d}$ 
 and $(\alpha,\beta,p)\in A\times B\times P$
the monotonicity condition
\begin{equation}
                                                     \label{5.23.1}
2\langle x-y,b^{\alpha\beta}(p,x)-b^{\alpha\beta}(p,y)\rangle
+\|\sigma^{\alpha\beta}(p,x)-
\sigma^{\alpha\beta}(p,y)\|^{2}\leq K_{1}|x-y|^{2},
\end{equation}
 holds,
where     $K_{1}$ is
a fixed constant. For instance, if
$\sigma^{\alpha\beta}(p,x)$ and $b^{\alpha\beta}(p,x)$
are Lipschitz continuous in $x$ with constant
independent of $\alpha,\beta,p$,
then \eqref{5.23.1} holds. If $d=1$, then \eqref{5.23.1}
is satisfied if, for instance $b^{\alpha\beta}(p,x)$ is
a decreasing function
and $\sigma^{\alpha\beta}(p,x)$  
is Lipschitz continuous in $x$ with constant
independent of $\alpha,\beta,p$.
Even if $\sigma$ and $b$ are independent of $p$,
this argument shows how 
control adapted processes may appear.

We discuss a different way in which
control adapted processes appear naturally
in Remark \ref{remark 5.6.1}.

\end{remark}

Take a $\zeta\in C^{\infty}_{0}(\bR^{d})$ with unit
integral and for $\varepsilon>0$ introduce
$\zeta_{\varepsilon}(x)=\varepsilon^{-d}\zeta(x/\varepsilon)$.
For locally summable functions $u=u(x)$ on $\bR^{d}$ define
$$
u^{(\varepsilon)}(x)=u*\zeta_{\varepsilon}(x).
$$

\begin{assumption}
                                             \label{assumption 5.23.1}

(i) For any $x \in\bR^{d}$  
\begin{equation}
                                              \label{3.27.4}
\sup_{(\alpha,\beta \in A\times B }(
|\bar c^{\alpha\beta}|+|\bar f^{\alpha\beta}|)( x)<\infty
\end{equation}

(ii)   There exist a constant $\delta_{1}\in(0,1]$ and a  function
$r^{\alpha\beta}(p,x)$ defined on $A\times B\times 
P\times \bR^{d}$ with values in $[\delta_{1},\delta_{1}^{-1}]$
such that $r^{\alpha\beta}(\bar{p},x)\equiv1$ and on
$A\times B\times 
P\times \bR^{d}$ we have
$$
  f ^{\alpha\beta}(p,x)=r^{\alpha\beta}(p,x)
 \bar f  ^{\alpha\beta}(x).
$$

(iii) For any bounded domain $D\subset\bR^{d}$ we have
$$
\|\sup_{(\alpha,\beta )\in A\times B }|\bar f^{\alpha\beta}  
 |\,\|_{L_{d}(D)}
+\|\sup_{(\alpha,\beta )\in A\times B } \bar c^{\alpha\beta } 
\,\|_{L_{d}(D)}<\infty,
$$
$$
\|\sup_{(\alpha,\beta )\in A\times B }|\bar f^{\alpha\beta} -
(\bar f^{\alpha\beta})^{(\varepsilon)}|\,\|_{L_{d}(D)}\to0,
$$
$$
 \|\sup_{(\alpha,\beta )\in A\times B }|\bar c^{\alpha\beta }
-(\bar c^{\alpha\beta})^{(\varepsilon)}|\,\|_{L_{d}(D)}\to0,
$$
  as 
$\varepsilon\downarrow0$.

(iv) There is a constant $\delta\in(0,1]$
such that for $\alpha\in A$, $\beta\in B$, $p\in P$,
and $x,\lambda\in\bR^{d}$ we have
$$
\delta|\lambda|^{2}\leq a^{\alpha\beta}_{ij}(p,x)\lambda_{i}
\lambda_{j}\leq \delta^{-1}|\lambda|^{2}.
$$

The reader understands, of course, that the summation
convention is adopted throughout the article.

\end{assumption}

Set
$$
\phi^{\alpha_{\cdot}\beta_{\cdot} x}_{t}
=\int_{0}^{t}c^{\alpha_{s}
\beta_{s} }( p^{\alpha_{\cdot}\beta_{\cdot}}_{s},x^{\alpha_{\cdot} 
\beta_{\cdot}  x}_{s})\,ds,
$$
fix a bounded domain $D\subset\bR^{d}$,
 define $\tau^{\alpha_{\cdot}\beta_{\cdot} x}$ as the first exit
time of $x^{\alpha_{\cdot} 
\beta_{\cdot} x}_{t}$ from $D$, and introduce
\begin{equation}
                                                    \label{2.12.2}
v(x)=\infsup_{\bbeta\in\bB\,\,\alpha_{\cdot}\in\frA}
E_{x}^{\alpha_{\cdot}\bbeta(\alpha_{\cdot})}\big[\int_{0}^{\tau}
f(p_{t},x_{t})e^{-\phi_{t}}\,dt+g(x_{\tau})e^{-\phi_{\tau}}\big],
\end{equation}
where the indices $\alpha_{\cdot}$, $\bbeta$, and $x$
at the expectation sign are written  to mean that
they should be placed inside the expectation sign
wherever and as appropriate, that is
$$
E_{x}^{\alpha_{\cdot}\beta_{\cdot}}\big[\int_{0}^{\tau}
f(p_{t},x_{t})e^{-\phi_{t}}\,dt+g(x_{\tau})e^{-\phi_{\tau}}\big]
$$
$$
:=
E \big[
 g(x^{\alpha_{\cdot}\beta_{\cdot}  x}
_{\tau^{\alpha_{\cdot}\beta_{\cdot}  x}}
)
e^{-\phi^{\alpha_{\cdot}\beta_{\cdot}  x}
_{\tau^{\alpha_{\cdot}\beta_{\cdot}  x}}}
+\int_{0}^{\tau^{\alpha_{\cdot}\beta_{\cdot}  x}}
f^{\alpha_{t}\beta_{t}   }
( 
p^{\alpha_{\cdot}\beta_{\cdot}}_{t},
x^{\alpha_{\cdot}\beta_{\cdot}  x}_{t})
e^{-\phi^{\alpha_{\cdot}\beta_{\cdot}  x}_{t}}\,dt\big].
$$
Observe that $v(x)=g(x)$ in $\bR^{d}\setminus D$.

This definition makes perfect sense due to the following.
\begin{lemma}
                                        \label{lemma 1.25.1}
There is a constant $N$, depending only
on $K_{0}$, $\delta$, $d$, and the diameter of $D$, such that for any
$\alpha_{\cdot}\in\frA$, $\beta_{\cdot}\in\frB$, $x\in D$, 
$n=1,2,...$, $t\in[0,\infty)$, and 
$h\in L_{d}(D)$ we have (a.s.)
\begin{equation}
                                                      \label{2.17.3}
I_{\tau^{\alpha_{\cdot}\beta_{\cdot} x}>t}E^{\alpha_{\cdot}\beta_{\cdot} }_{x}\big\{
\big(\int_{t}^{\tau}|h(x_{s})|
\,ds\big)^{n}\mid\cF_{t}\big\}\leq n!N^{n}\|h\|_{L_{d}(D)}^{n}.
\end{equation}
In particular, for any
$n=1,2,...$
$$
E^{\alpha_{\cdot}\beta_{\cdot} }_{x}\tau^{n}\leq n! N^{n}.
$$

\end{lemma}

 Proof. Estimate \eqref{2.17.3} with $n=1$ is proved
in Theorem 2.2.1  of \cite{Kr77}. If it is true for an $n$,
then we have
$$
I_{\tau^{\alpha_{\cdot}\beta_{\cdot} x}>t}
E^{\alpha_{\cdot}\beta_{\cdot} }_{x}\big\{
\big(\int_{t}^{\tau}|h(x_{s})|
\,ds\big)^{n+1}\mid\cF_{t}\big\}
$$
$$
=(n+1)I_{\tau^{\alpha_{\cdot}\beta_{\cdot} x}>t}
E^{\alpha_{\cdot}\beta_{\cdot} }_{x}\big\{
\int_{t}^{\tau}|h(x_{r})|\big(\int_{r}^{\tau}|h(x_{s})|
\,ds\big)^{n}\,dr \mid\cF_{t}\big\}
$$
$$
=(n+1)I_{\tau^{\alpha_{\cdot}\beta_{\cdot} x}>t}
E^{\alpha_{\cdot}\beta_{\cdot} }_{x}\big\{
\int_{t}^{\tau}|h(x_{r})|I_{\tau>r}
\big[E^{\alpha_{\cdot}\beta_{\cdot} }_{x}\big(\int_{r}^{\tau}|h(x_{s})|
\,ds\big)^{n} \mid\cF_{r}\big]\,dr  \mid\cF_{t}\big\}
$$
$$
\leq N^{n}(n+1)!\|h\|^{n}_{L_{d}(D)}
I_{\tau^{\alpha_{\cdot}\beta_{\cdot} x}>t}
E^{\alpha_{\cdot}\beta_{\cdot} }_{x}\big\{
\int_{t}^{\tau}|h(x_{r})|\,dr  \mid\cF_{t}\big\}
$$
$$
\leq  N^{n+1}(n+1)!\|h\|^{n+1}_{L_{d}(D)}.
$$
The lemma is proved.

For a sufficiently smooth function $u=u(x)$ introduce
$$
L^{\alpha\beta} u(p,x)=a^{\alpha\beta}_{ij}(p,x)D_{ij}u(x)+
b ^{\alpha\beta}_{i }(p,x)D_{i}u(x)-c^{\alpha\beta} (p,x)u(x),
$$
where, naturally, $D_{i}=\partial/\partial x_{i}$, $D_{ij}=D_{i}D_{j}$.
Recall that we fixed a $\bar{p}\in P$ and denote
$$
\bar{L}^{\alpha\beta} u(x)=L^{\alpha\beta}u( \bar{p},x),
$$
\begin{equation}
                                                     \label{1.16.1}
H[u](x)=\supinf_{\alpha\in A\,\,\beta\in B}
[\bar{L}^{\alpha\beta} u(x)+\bar{f}^{\alpha\beta} (x)].
\end{equation}

\begin{definition}
                                        \label{definition 5.8.1}

 For a domain  $U\subset \bR^{d}$ we say that a
$C^{2}_{loc}(U)$ function $u$ is
$p$-insensitive in $U$ (relative to
$(r^{\alpha\beta}, L^{\alpha\beta}) $)
if  
for any $x\in U$, $\alpha_{\cdot}\in\frA$, and $\beta_{\cdot}
\in\frB$  
$$
d\big[u(x^{\alpha_{\cdot}\beta_{\cdot}x}_{t})
e^{-\phi^{\alpha_{\cdot}\beta_{\cdot}  x}_{t}} \big]=
r^{\alpha_{t}\beta_{t}}(p_{t}
^{\alpha_{\cdot}\beta_{\cdot}},
x^{\alpha_{\cdot}\beta_{\cdot}  x}_{t})
  \bar{L}^{\alpha_{t}\beta_{t}}
u(x^{\alpha_{\cdot}\beta_{\cdot}x}_{t})
e^{-\phi^{\alpha_{\cdot}\beta_{\cdot}  x}_{t}}\,dt +dm_{t}
$$
for $t$ less than the first exit time of 
$x^{\alpha_{\cdot}\beta_{\cdot}  x}_{t}$ from $U$, where
$m_{t}$ is a local martingale starting at zero.

\end{definition}

There are nontrivial cases when all sufficiently smooth
functions are $p$-insensitive (see Example \ref{example 5.8.1}).
On the other hand, any smooth function $u(x_{1})$ will be 
$p$-insensitive if $(a_{11},b_{1})^{\alpha\beta}(p,x)=
r ^{\alpha\beta}(p,x)
(\bar{a}_{11},\bar{b}_{1})^{\alpha\beta}(x)$ with {\em no\/}
restrictions on other entries of $a$ and $b$. A generalization
of this particular example will play an extremely important
role in one of subsequent articles.

\begin{definition}
Let $U$ be a domain in $\bR^{d}$ for which the Sobolev
embedding $W^{2}_{d}(U)\subset C(\bar U)$ is valid.
 We say that it is regular
(for given $g$)
if there exists a    function
$u\in W^{2}_{d}(U)$ such that $H[u]=0$ in $U$ (a.e.) and
$u=g$ on $\partial U$ and there exists a sequence 
$u_{n}\in C^{2}(\bar{U})$ of  $p$-insensitive in $U$ functions
such that $u_{n}\to u$
in $W^{2}_{d}(U)$ and in $C(\bar{U})$.

\end{definition}

In a subsequent article we will show that the following assumption
can be dropped.

\begin{assumption}
                                           \label{assumption 3.22.1}
There is a sequence of expanding regular subdomains
$D_{n}$ of $D$ such that $D=\bigcup_{n\geq1}D_{n}$.

\end{assumption}

Finally we impose the following.

\begin{assumption}
                                           \label{assumption 3.19.1}
There exists a bounded
 nonnegative $G\in   C^{2}_{loc}(D)$
such that 

(i) We have $G\in C(\bar{D})$ and $G=  0$ on $\partial D$;

(ii) For all $\alpha\in A$, $\beta\in B$, $p\in P$,
and $x\in D$
\begin{equation}
                                                    \label{3.20.1}
  L^{\alpha\beta}G(p,x)\leq-1.
\end{equation}

\end{assumption}

Here is our main result.

\begin{theorem}
                                                    \label{theorem 1.14.1}
 
Under the above assumptions also suppose
that there exists a sequence  of
functions $g_{n}$ such that $\|g-g_{n}\|_{C(\bar{D})}\to0$
as $n\to \infty$, for each $n\geq1$,
$\|g_{n}\|_{ C^{2}(\bar{D})}<\infty$ and $g_{n}$ is
$p$-insensitive in $D$.
 Then

(i) The function $v(x)$  
is independent
of the chosen control adapted process $p\in\frP$,
it is bounded and continuous in $\bR^{d}$.

(ii)
Let   $\gamma^{\alpha_{\cdot}\beta_{\cdot}x} $
 be an $\{\cF_{t}\}$-stopping
time defined for each $\alpha_{\cdot}\in\frA$, $\beta\in\frB$,
and $ x\in\bR^{d}$
 and such
 that $\gamma^{\alpha_{\cdot}\beta_{\cdot}x}\leq 
\tau^{\alpha_{\cdot}\beta_{\cdot} x}$. Also let 
$\lambda_{t}^{\alpha_{\cdot}
\beta_{\cdot}x}\geq0$ be progressively measurable  functions 
on $\Omega
\times[0,\infty)$ 
defined for each $\alpha_{\cdot}\in\frA$, $\beta\in\frB$,
and $ x\in\bR^{d}$ and
such that they have finite integrals over finite time intervals
(for any $\omega$).
Then for any $x$
\begin{equation}
                                                   \label{1.14.1}
v(x)=\infsup_{\bbeta\in\bB\,\,\alpha_{\cdot}\in\frA}
E_{x}^{\alpha_{\cdot}\bbeta(\alpha_{\cdot})}\big[
v(x_{\gamma})e^{-\phi_{\gamma}-\psi_{\gamma}}
+\int_{0}^{\gamma}
\{f(p_{t},x_{t})+\lambda_{t}v(x_{t})\}e^{-\phi_{t}-\psi_{t}}\,dt \big],
\end{equation}
where inside the expectation sign
$\gamma=\gamma^{\alpha_{\cdot}\bbeta(\alpha_{\cdot})x}
$ and
$$
\psi^{\alpha_{\cdot}\beta_{\cdot}x} _{t}
=\int_{0}^{t}
\lambda^{\alpha_{\cdot}\beta_{\cdot}x}_{s}\,ds.
$$

\end{theorem}

\begin{remark}
The function $G$ is called a barrier in the theory 
of partial differential equations. Existence
of such barriers is known for a very large class of domains,
say such that there are  $\rho_{0}>0$ and $\theta>0$ such that
for any point $x_{0}\in\partial D$ and any $r\in(0,\rho_{0}]$
we have that the volume of the intersection of $D^{c}$
with the ball of radius $r$ centered at $x_{0}$ is greater
than $\theta r^{d}$. The so-called uniform exterior cone condition
will suffice.

Without Assumption \ref{assumption 3.19.1}
or similar ones one cannot assert that
$v$ is continuous
 in $\bar{D}$ even if no control parameters are involved.
\end{remark}

Note that
the possibility to vary $\lambda$ in Theorem \ref{theorem 1.14.1}
might be useful while considering stochastic  differential  games
with stopping in the spirit of \cite{Kr70}.

\begin{remark} 
                                 \label{remark 5.11.1}
Definition \ref{definition 5.8.1} is stated in the form
which is easy to use and to check especially
when (as in a subsequent article) the state
process consists of several components for each of which
the corresponding equations have very different
forms and $u$ depends only on part of these components.

Still it is worth noting that, as follows immediately from
It\^o's formula, $u\in C^{2}_{loc}(U)$
 will be $p$-insensitive if
on $A\times B\times 
P\times U$ we have  
\begin{equation}
                                                  \label{4.11.2}
L^{\alpha\beta} u(p,x) =
r^{\alpha\beta}(p,x) \bar{L}^{\alpha\beta}  u( x).
\end{equation}

\end{remark}
\begin{example}
                                                 \label{example 5.8.1}
Let $\cO$ be the set  
of $d_{1}\times d_{1}$ orthogonal matrices and denote by $p=(p',p'')$
a generic point in $[-\delta_{1},\delta_{1}^{-1}]\times\cO$. 
Assume that the original $\sigma,b$ are independent of $p$.
Then introduce
$$
\sigma^{\alpha\beta}(p,x)=(p')^{1/2}\sigma^{\alpha\beta}(x)p'',\quad
(b,c,f)^{\alpha\beta}(p,x)=p'(b,c,f)^{\alpha\beta}(x).
$$
In this case
$$
(1/2)\sigma^{\alpha\beta}(p,x)(\sigma^{\alpha\beta}(p,x))^{*}
=p'a^{\alpha\beta}(x)
$$
and \eqref{4.11.2} holds for any $u$ with $r^{\alpha\beta}(p,x)=p'$.

In this case the assertion that $v$
is independent
of the control adapted process $p_{t}^{\alpha_{\cdot}\beta_{\cdot}}$
is rather natural and is due to the fact that
its effect on the state process is equivalent to that
of  a random time change
and a random rotation of the increments of the original Wiener process.

The main advantage of  introducing the above parameters, which
by far are not the most general and important for the future,
is that while  
 estimating $v(x+\varepsilon\xi)-v(x)$ for small $\varepsilon$,
 where $\xi\in \bR^{d}$, we can take $p\equiv (1,I)$, where $I$
is the $d_{1}\times d_{1}$ identity matrix, in the definition
of $v(x)$ and a different $p$ close to $(1,I)$ in the definition
of $v(x+\varepsilon\xi)$. This may make the solutions
of the corresponding stochastic equations to become closer 
than in  
the case where for both $v(x )$
and $v(x+\varepsilon\xi)$ we   take $p\equiv (1,I)$.

\end{example}
 
\begin{remark}
                                           \label{remark 5.6.1}
One of ways to introduce control adapted processes
can be explained in the situation of Example \ref{example 5.8.1}
when the original $\sigma$ and $b$ are Lipschitz continuous.
Take a $[\delta_{1},\delta_{1}^{-1}]$-valued
function $r(x)$ and $\cO$-valued function $Q(x)$
defined and Lipschitz continuous on $\bR^{d}$.
Fix an $x_{0}\in\bR^{d}$ and for
 $\alpha_{\cdot}\in\frA$ and $\beta_{\cdot}\in
\frB$ define
$$
p^{\alpha_{\cdot}\beta_{\cdot}}_{t}:=
(r(x_{t}),Q(x_{t})),
$$
where $x_{t}$ is a unique solution of
$$
x_{t}=x_{0}+\int_{0}^{t}r^{1/2}(x_{s})
\sigma^{\alpha_{s}\beta_{s}}(x_{s})Q(x_{s})\,dw_{s}
+\int_{0}^{t}r(x_{s})
b^{\alpha_{s}\beta_{s}}(x_{s})\,ds.
$$
Almost obviously $p\in\frP$ and the above solution
is, actually, $x^{\alpha_{\cdot}\beta_{\cdot}x}_{t}$
for that $p$ if $x=x_{0}$. In a subsequent article we will show
a much more sophisticated use of control adapted processes
defined by  an auxiliary It\^o  equation.
\end{remark}

As a simple byproduct of our proofs we obtain the following.

\begin{theorem}
                                             \label{theorem 2.18.1}
The function $v$ is locally H\"older continuous
in $D$ with exponent $\theta\in(0,1)$ depending only
on $d$ and $\delta$.

\end{theorem}

The point is that $v$ will be obtained as the limit of $u_{n}$
which are solutions of class $W^{2}_{d}(D_{n})$
of the equation $H[u_{n}]=0$ in $D_{n}$ (a.e.) 
with boundary data $g$.
It is well known (see, for instance, Remark 1.3 in
\cite{Kr12.2})
 that such $u_{n}$ 
satisfy linear uniformly elliptic equations
with bounded coefficients and it  is a classical result that
such solutions admit uniform in $n$ local H\"older
estimates of some exponent $\theta\in(0,1)$ depending only
on $d$ and $\delta$ (see, for instance,
 \cite{GT} or \cite{Kr85}).

\mysection{Proof of Theorem \protect\ref{theorem 1.14.1}
in case that the Isaacs equation has a smooth solution}

                                          \label{section 2.11.1}

In this section Assumptions \ref{assumption 5.23.1} (iii), (iv),  
 the regularity Assumption
\ref{assumption 3.22.1} as well as Assumption
\ref{assumption 3.19.1} concerning  $G$ 
are not used  and the domain $D$ 
  is not supposed to be bounded.
Suppose that for each $\varepsilon>0$ we are given real-valued
functions $c^{\alpha\beta}_{\varepsilon}(x )$ and
$f^{\alpha\beta}_{\varepsilon}(x )$ defined on $A\times B\times
\bR^{d} $. 
\begin{assumption}
                                           \label{assumption 2.19.1}
(i) Assumptions \ref{assumption 1.9.1}, \ref{assumption 5.23.2},   
and \ref{assumption 5.23.1} (i), (ii)
are satisfied.

(ii) For a constant $\chi>0$ we have
$c^{\alpha\beta}(p,x )\geq\chi$ for all $\alpha,\beta,p,x$.

(iii) For each $\varepsilon>0$   
the functions $(c,f)_{\varepsilon}^{\alpha\beta}(x )$ are 
bounded on $A\times B\times \bar{D}$ and uniformly
 continuous with respect to $x\in\bar{D}$
uniformly with respect to $\alpha,\beta $.

(iv) For any $x$ as $\varepsilon\downarrow0$,  
$$
d_{\varepsilon}(x):
=\sup_{(\alpha_{\cdot},\beta_{\cdot} )\in\frA\times\frB }
E^{\alpha_{\cdot}\beta_{\cdot}}_{x}
\int_{0}^{\tau}(|\bar{c}-c_{\varepsilon}|+|
\bar{f}-f_{\varepsilon}|)(x_{t} )
e^{-\phi_{t}}\,dt\to0.
$$

(v) For any $x\in D$
$$
\sup_{(\alpha_{\cdot}\beta_{\cdot})\in\frA\times\frB}
E^{\alpha_{\cdot}\beta_{\cdot}}_{x}
\int_{0}^{\tau}|f (p_{t},x_{t})|e^{-\phi_{t}}\,dt<\infty.
$$

\end{assumption}
 
Observe that Assumption \ref{assumption 2.19.1} (v) implies that
$v$ is well defined.
 
In some applications we have in mind the following
``degenerate" version of Theorem \ref{theorem 1.14.1} plays
an important role.
We assume that we are given two $p$-insensitive in $D$ functions
 $\hat{u},\check{u} \in C^{2}(\bar{D})$
(with finite $ C^{2}(\bar{D})$-norms) such that their
second-order derivatives are   uniformly continuous in $\bar D$
(in case $D$ is unbounded). 

\begin{theorem}
                                         \label{theorem 2.11.1}

(i) If $H[\hat{u}]\leq0$ (everywhere)
  in $D$ 
and $\hat{u}\geq g$ on $\partial D$
(in case $\partial D\ne\emptyset$), then $v\leq\hat{u}$ in $\bar{D}$.

(ii) If $H[\check{u}]\geq0$ (everywhere)
  in $D$ 
and $\check{u}\leq g$  on $\partial D$
(in case $\partial D\ne\emptyset$), then $v\geq\check{u}$ in $\bar{D}$.

(iii) If $\hat{u}$ and $\check{u}$ are as in (i) and (ii) and
  $\hat{u}=\check{u}$, then
all assertions of Theorem \ref{theorem 1.14.1} hold true. Moreover,
  $v=\hat{u}$.
\end{theorem}

This theorem is an immediate consequence of the following two results
in which we allow $\lambda^{\alpha_{\cdot}\beta_{\cdot}x}_{t}$
and $\gamma^{\alpha_{\cdot}\beta_{\cdot}x}$
to be as in Theorem \ref{theorem 1.14.1}.
\begin{theorem}
                                         \label{theorem 4.17.1} 
Suppose that $H[\hat{u}]\leq0$ (everywhere) in $D$. Then
for all $x\in\bar{D}$ we have
$$
\hat{u}(x)\geq \infsup_{\bbeta\in\bB\,\,\alpha_{\cdot}\in\frA}
E_{x}^{\alpha_{\cdot}\bbeta (\alpha_{\cdot})}
\big[\hat{u}(x_{\gamma})e^{-\phi_{\gamma}-\psi_{\gamma}}
$$
\begin{equation}
                                                     \label{2.11.5}
+
 \int_{0}^{\gamma}
[f(x_{t},p_{t})+\lambda_{t}\hat{u}(x_{t}) ]e^{-\phi_{t}-\psi_{t}}\,dt\big].
\end{equation}
In particular,  
if $\hat{u}\geq g$ on $\partial D$, then
for $\gamma\equiv\tau$ and $\lambda\equiv0$
equation \eqref{2.11.5} yields that $\hat{u}\geq v$.
 
\end{theorem}
\begin{theorem}
                                         \label{theorem 4.17.2} 
Suppose that $H[\check{u}]\geq0$ (everywhere)
  in $D$. Then
for all $x\in\bar{D}$ we have
$$
\check{u}(x)\leq \infsup_{\bbeta\in\bB\,\,\alpha_{\cdot}\in\frA}
E_{x}^{\alpha_{\cdot}\bbeta (\alpha_{\cdot})}
\big[\check{u}(x_{\gamma})e^{-\phi_{\gamma}-\psi_{\gamma}}
$$
\begin{equation}
                                                     \label{2.11.05}
+
 \int_{0}^{\gamma}
[f(p_{t},x_{t})+\lambda_{t}\check{u}(x_{t}) ]e^{-\phi_{t}-\psi_{t}}\,dt\big].
\end{equation}
In particular, if  $\check{u}\leq g$ on $\partial D$,
then for $\gamma\equiv\tau$ and $\lambda\equiv0$
equation \eqref{2.11.05} yields that $\check{u}\leq v$. 
 
\end{theorem}

Note that, formally, the value $x_{\gamma}$ 
in \eqref{2.11.5} and \eqref{2.11.05} 
 may not be defined if $\gamma=\infty$.
In that case we set the corresponding terms to equal zero,
which is natural because $\hat{u}$ and $\check{u}$
are bounded and $\phi^{\alpha_{\cdot}\beta_{\cdot}x}_{\infty}
=\infty$.

To prove these theorems we need two lemmas.
The reader can compare our arguments with the ones
in \cite{Sw96} and see that they are very close.

For a stopping time $\gamma$
we say that a process $\xi_{t}$ is a submartingale
on $[0,\gamma]$ if $\xi_{t\wedge\gamma}$ is a submartingale.
Similar definition applies to supermartingales.

\begin{lemma}
                                                  \label{lemma 1.16.1}
Let $H[\hat{u}]\leq0$ (everywhere)
  in $D$. Then for any $x\in\bR^{d}$,
 $\alpha_{\cdot}\in\frA$, and $\varepsilon>0$,
there exist a sequence $\beta^{n}_{\cdot}(\alpha_{\cdot})=
\beta^{n}_{\cdot}(\alpha_{\cdot},x,\varepsilon)
\in\frB$, $n=1,2,...$, and a sequence of increasing continuous
$\{\cF_{t}\}$-adapted processes 
$\eta^{n\varepsilon}_{t}(\alpha_{\cdot})
=\eta^{n\varepsilon}_{t}(\alpha_{\cdot},x)$
 with $\eta^{n\varepsilon}_{0}(\alpha_{\cdot})=0$
 such that 
\begin{equation}
                                                 \label{2.29.1}
\sup_{n}E 
\eta^{n\varepsilon}_{\infty}(\alpha_{\cdot})<\infty,
\end{equation}
 the processes
$$
\kappa^{n\varepsilon}_{t}(\alpha_{\cdot}):=
\hat{u}(x^{n}_{t})
e^{-\phi^{n}_{t}}
-\eta^{n\varepsilon}_{t} (\alpha_{\cdot})
+\int_{0}^{t}
f^{n}_{s} ( p^{n}_{s},x^{n}_{s}) 
e^{-\phi^{n}_{s}}\,ds,
$$
where
\begin{equation}
                                                   \label{5.6.1}
(x^{n}_{t},\phi^{n}_{t})=(x_{t},\phi_{t})
^{\alpha_{\cdot}\beta^{n}_{\cdot}(\alpha_{\cdot})  x},\quad
f^{n}_{t}(p,x)=f^{\alpha_{t}\beta^{n}_{t}(\alpha_{\cdot})}(p,x),\quad
p^{n}_{t}=p^{\alpha_{\cdot}\beta^{n}_{\cdot}(\alpha_{\cdot})}_{t}.
\end{equation}
are supermartingales on $[0,\tau^{\alpha_{\cdot}\beta^{n}_{\cdot}(\alpha_{\cdot})
  x}]$, and
\begin{equation}
                                                       \label{3.3.1}
\nlimsup_{n\to\infty}
\sup_{\alpha_{\cdot}\in\frA}E 
 \eta^{n\varepsilon}_{\tau}(\alpha_{\cdot})\leq
  \varepsilon/(\delta_{1}\chi) +Nd_{\varepsilon}(x),
\end{equation}
where $\delta_{1}$
is taken from Assumption \ref{assumption 5.23.1} (ii) and
$N$ is independent of $x$ and $\varepsilon$.
 Furthermore, if for any $n$ we are given a nonnegative,
 progressively measurable process 
$\lambda^{n}_{t}\geq0$ having finite integrals over
finite time intervals (for any $\omega$),
then the processes  
$$
\rho^{n\varepsilon}_{t}(\alpha_{\cdot}):=
\hat{u}(x^{n}_{t})
e^{-\phi^{n}_{t}
-\psi^{n}_{t}}
-\eta^{n\varepsilon}_{t} (\alpha_{\cdot})e^{-\psi^{n}_{t}}
$$
$$
+\int_{0}^{t}
\big[f^{n} _{s}
(p^{n}_{s}
, x^{n}_{s})
+\lambda^{n}_{s}\hat{u}
(x^{n}_{s})
-\lambda^{n}_{s}\eta^{n\varepsilon}_{s} (\alpha_{\cdot})
e^{ \phi^{n}_{s}}\big]
e^{-\phi^{n}_{s}
-\psi^{n}_{t}}\,ds 
$$
are supermartingales on 
$[0,\tau^{\alpha_{\cdot}\beta^{n}_{\cdot}(\alpha_{\cdot})
  x}]$, where (we use  notation \eqref{5.6.1} and)
\begin{equation}
                                                    \label{4.28.1}
\psi^{n}_{t}=\int_{0}^{t}\lambda^{n}_{s}\,ds.
\end{equation}
 Finally,
\begin{equation}
                                           \label{2.29.3}
\sup_{\alpha_{\cdot}\in\frA}\sup_{n}
E 
\sup_{t\geq0}|
\kappa^{n\varepsilon}_{t\wedge\tau }(\alpha_{\cdot})|<\infty,\quad
 \sup_{\alpha_{\cdot}\in\frA}\sup_{n}
E  
\sup_{t\geq0}|
\rho^{n\varepsilon}_{t\wedge\tau }(\alpha_{\cdot})|<\infty.
\end{equation}

\end{lemma}  

Proof. Since $B$ is separable and $a^{\alpha\beta},b^{\alpha\beta},
c^{\alpha\beta}$, and $f^{\alpha\beta}$ are continuous with respect to
$\beta$ one can replace $B$ in \eqref{1.16.1} with an appropriate
countable subset $B_{0}=\{\beta_{1},\beta_{2},...\}$.
 Then for each $\alpha\in A$ and $x\in D$
define $\beta(\alpha,x)$ as $\beta_{i}\in B_{0}$ with the least 
$i$ such that
\begin{equation}
                                                  \label{1.17.2}
0
\geq \bar{L}^{\alpha\beta_{i}} \hat{u}(x)+\bar{f}^{\alpha\beta_{i}} 
(x)-\varepsilon.  
\end{equation}
For each $i$ the right-hand side of \eqref{1.17.2} is
Borel in $x$ and continuous in $\alpha$. Therefore,
it is a Borel function of $(\alpha,x)$, implying that  $\beta(\alpha,x)$ 
also is a Borel function
of $(\alpha,x)$. For $x\not\in D$ set $\beta(\alpha,x)=\beta^{*}$,
where $\beta^{*}$ is a fixed element of $B$.
Then we have that in $D$
\begin{equation}
                                                  \label{1.17.3}
0
\geq \bar{L}^{\alpha\beta(\alpha,x)} \hat{u}(x)
+\bar{f}^{\alpha\beta(\alpha,x)} (x)-\varepsilon.
\end{equation}

 After that fix $x$,
define  $\beta^{n0}_{t}(\alpha_{\cdot})=\beta(\alpha_{t},x)$, $t\geq0$,
and for $k\geq1$ introduce
$\beta^{nk}_{t}(\alpha_{\cdot})$ 
recursively so that
\begin{equation}
                                                  \label{4.5.2}
\beta^{nk}_{t}(\alpha_{\cdot})=\beta^{n(k-1)}_{t}(\alpha_{\cdot})
\quad\text{for}\quad t<k/n,
\end{equation}
$$
\beta^{nk}_{t}(\alpha_{\cdot})=\beta(\alpha_{t},x^{nk }_{k/n})
\quad\text{for}\quad t\geq k/n,
$$
where $x^{nk}_{t}$ is a unique solution of
$$
x _{t}=x+\int_{0}^{t}\sigma(\alpha_{s},\beta^{n(k-1)}_{s}(\alpha_{\cdot}),
p^{\alpha_{\cdot}\beta^{n(k-1)}_{\cdot}(\alpha_{\cdot})}_{s},x_{s})\,
dw_{s}
$$
\begin{equation}
                                                      \label{4.5.1}
+\int_{0}^{t}b(\alpha_{s},\beta^{n(k-1)}_{s}(\alpha_{\cdot}),
p^{\alpha_{\cdot}\beta^{n(k-1)}_{\cdot}(\alpha_{\cdot})}_{s},x_{s})\,
ds.
\end{equation}

To show that the above definitions make sense, observe that, 
by Assumption \ref{assumption 5.23.2},
$x^{n0}_{t}$ is well defined for all $t$.
 Therefore, $\beta^{n1}_{t}(\alpha_{\cdot})$
is also well defined, and by induction we conclude that
$x^{nk}_{t}$ and $\beta^{nk}_{t}(\alpha_{\cdot}) $ are
well defined for all $k$.

Furthermore, owing to \eqref{4.5.2} it makes sense to define
$$
\beta^{n}_{t}(\alpha_{\cdot})=\beta^{nk}_{t}(\alpha_{\cdot})
\quad\text{for}\quad t<k/n.
$$
Notice that by definition 
 $x^{n}_{t}=x^{\alpha_{\cdot}\beta^{n}_{\cdot}(\alpha_{\cdot})x}_{t}$
satisfies the equation
$$
x _{t}=x+\int_{0}^{t}\sigma(\alpha_{s},\beta^{n }_{s}(\alpha_{\cdot}),
p^{\alpha_{\cdot}\beta^{n }_{\cdot}(\alpha_{\cdot})}_{s},x_{s})\,
dw_{s}
$$
\begin{equation}
                                                      \label{4.6.1}
+\int_{0}^{t}b(\alpha_{s},\beta^{n }_{s}(\alpha_{\cdot}),
p^{\alpha_{\cdot}\beta^{n }_{\cdot}(\alpha_{\cdot})}_{s},x_{s})\,
ds.
\end{equation}
For $t<k/n$ we have $\beta^{n}_{t}(\alpha_{\cdot})
=\beta^{n(k-1)}_{t}(\alpha_{\cdot})$,
so that for $t\leq k/n$ equation \eqref{4.6.1}
coincides with \eqref{4.5.1} owing to the
fact that  $p^{\alpha_{\cdot}
\beta_{\cdot}}_{t}$ is control adapted.
 It follows that (a.s.)
$$
x_{t}^{n}=
x^{nk}_{t}\quad\text{for all}\quad t\leq k/n,
$$
so that   (a.s.)
$$
\beta^{nk}_{t}(\alpha_{\cdot})=\beta(\alpha_{t},x^{n}_{k/n})
$$
for all $t\geq k/n$. Therefore, if $(k-1)/n\leq t<k/n$, then
$$
\beta^{n}_{t}(\alpha_{\cdot})=\beta^{n(k-1)}_{t}
(\alpha_{\cdot})=\beta(\alpha_{t},x^{n}_{(k-1)/n})
$$

\begin{equation}
                                                         \label{1.19.1}
\beta^{n}_{s}:=\beta^{n}_{s}(\alpha_{\cdot})=\beta(\alpha_{s},
x^{n}_{\kappa_{n}(s)}),
\end{equation}
where $\kappa_{n}(t)=[nt]/n$, and $x_{t}^{n}$
satisfies 
$$
x_{t}^{n}=x+\int_{0}^{t}\sigma(\alpha_{s},\beta(\alpha_{s},
x^{n}_{\kappa_{n}(s)}),p^{n}_{s},x^{n}_{s})\,dw_{s}
$$
\begin{equation}
                                           \label{1.17.1}
+\int_{0}^{t}b(\alpha_{s},\beta(\alpha_{s},
x^{n}_{\kappa_{n}(s)}),p^{n}_{s},x^{n}_{s})\,ds,
\end{equation}
with $p^{n}_{s}
=p^{\alpha_{\cdot}\beta^{n }_{\cdot} }_{s}$.  

 Introduce $\tau^{n}$ as the first exit time  
of $x^{n}_{t}$ from $D$ and set
$$
\phi^{n}_{t}=\phi^{\alpha_{\cdot}\beta^{n}_{\cdot}x}_{t},\quad
r^{n}_{s}=r^{\alpha_{s}\beta^{n}_{s}}(p^{n}_{s},x^{n}_{s}).
$$

 Observe that
by It\^o's formula
$$
\hat{u}(x^{n}_{t\wedge\tau^{n}})e^{-\phi^{n}_{t\wedge \tau^{n}}}
=\hat{u}(x)+\int_{0}^{t\wedge\tau^{n}}e^{-\phi^{n}_{s}}
 L^{\alpha_{s}\beta^{n}_{s}}
\hat{u}(p^{n}_{s},x^{n}_{s})
\,ds+m^{n}_{t} ,
$$where $m^{n}_{s}$ is a martingale.

 By Definition \ref{definition 5.8.1}
$$
\hat{u}(x^{n}_{t\wedge\tau^{n}})e^{-\phi^{n}_{t\wedge \tau^{n}}}
+\int_{0}^{t\wedge \tau^{n}}f^{\alpha_{s}\beta^{n}_{s}}(p^{n}_{s},x^{n}_{s})
e^{-\phi^{n}_{s}}\,ds
$$
\begin{equation}
                                                         \label{5.8.5}
=\hat{u}(x)+\int_{0}^{t\wedge\tau^{n}}e^{-\phi^{n}_{s}}
 r^{n}_{s}\big[
 \bar{L}^{\alpha_{s}\beta^{n}_{s}}
\hat{u}( x^{n}_{s})+\bar{f}^{\alpha_{s}\beta^{n}_{s}}
 ( x^{n}_{s})\big]
\,ds+m^{n}_{t} ,
\end{equation}
where,
for $s<\tau^{n}$,
(notice the change of $\bar c$ to $c_{\varepsilon}$)
$$
\bar L^{\alpha_{s}\beta^{n}_{s}}
\hat{u}(x^{n}_{s})=\bar a_{ij}(\alpha_{s},\beta(\alpha_{s},
x^{n}_{\kappa_{n}(s)}),x^{n}_{s})
D_{ij}\hat{u}(x^{n}_{s})
$$
$$
+\bar b_{i}(\alpha_{s},\beta(\alpha_{s},x^{n}_{\kappa_{n}(s)}),x^{n}_{s})
D_{i}\hat{u}(x^{n}_{s})
-\bar c(\alpha_{s},\beta(\alpha_{s},x^{n}_{\kappa_{n}(s)}),x^{n}_{s})
 \hat{u}(x^{n}_{s})
$$
$$
=\bar a_{ij}(\alpha_{s},\beta(\alpha_{s},
x^{n}_{\kappa_{n}(s)}),x^{n}_{s})
D_{ij}\hat{u}(x^{n}_{s})+
\bar b_{i}(\alpha_{s},\beta(\alpha_{s},x^{n}_{\kappa_{n}(s)}),x^{n}_{s})
D_{i}\hat{u}(x^{n}_{s})
$$
$$
-c_{\varepsilon}(\alpha_{s},\beta(\alpha_{s},x^{n}_{\kappa_{n}(s)}),x^{n}_{s})
 \hat{u}(x^{n}_{s})+\xi^{n\varepsilon}_{t},
$$
where $\xi^{n\varepsilon}_{t}$ (defined by the above equality)
is a progressively measurable process such that
by Assumption \ref{assumption 2.19.1} (iv)
\begin{equation}
                                                       \label{2.29.2}
\sup_{n}E\int_{0}^{\tau^{n}}|\xi^{n\varepsilon}_{t}|e^{-\phi^{n}_{t}}
\,dt\leq Nd_{\varepsilon}(x)
\end{equation}
with $N$   independent of $\alpha$, $\varepsilon$, and $x$
(equal to one).  All such processes are denoted
by $\xi^{n\varepsilon}_{t}$ below even if they may change 
from one occurrence to another.

According to our assumptions on the uniform continuity in $x$ of the data 
and $D_{ij}u(x)$ we  have that   
for $s<\tau^{n}$ 
$$
\bar L^{\alpha_{s}\beta^{n}_{s}}
\hat{u}(x^{n}_{s})\leq  
\bar a_{ij}(\alpha_{s},\beta(\alpha_{s},
x^{n}_{\kappa_{n}(s)}),x^{n}_{\kappa_{n}(s)})
D_{ij}\hat{u}(x^{n}_{\kappa_{n}(s)})
$$
$$
+\bar b_{i}(\alpha_{s},\beta(\alpha_{s},x^{n}_{\kappa_{n}(s)}),x^{n}_{\kappa_{n}(s)})
D_{i}\hat{u}(x^{n}_{\kappa_{n}(s)})
$$
$$
-\bar c(\alpha_{s},\beta(\alpha_{s},x^{n}_{\kappa_{n}(s)}),x^{n}_{\kappa_{n}(s)})
 \hat{u}(x^{n}_{\kappa_{n}(s)}) 
+\chi_{\varepsilon}
(x^{n}_{s}-x^{n}_{\kappa_{n}(s)})+|\xi^{n\varepsilon}_{t}|
$$
where,
 for each $\varepsilon>0$, $\chi_{\varepsilon}(y)$ 
is a (nonrandom) bounded function on $\bR^{d}$ such that
$\chi_{\varepsilon}(y)\to0$ as $y\to0$. All 
such functions will be denoted
by $\chi_{\varepsilon}$ even if they may change from one occurrence to another.
Then  \eqref{1.17.3} shows that,
for $s<\tau^{n}$,
$$
\bar L^{\alpha_{s}\beta^{n}_{s}}
\hat{u}(x^{n}_{s})\leq \varepsilon+\chi_{\varepsilon}
(x^{n}_{s}-x^{n}_{\kappa_{n}(s)})
+|\xi^{n\varepsilon}_{t}|
 -\bar f(\alpha_{s},\beta(\alpha_{s}(x^{n}_{\kappa_{n}(s)}),
x^{n}_{\kappa_{n}(s)})
$$
$$
\leq \varepsilon+\chi_{\varepsilon}
(x^{n}_{s}-x^{n}_{\kappa_{n}(s)}) 
+|\xi^{n\varepsilon}_{t}| -\bar f^{\alpha_{s}\beta^{n}_{s}}
( x^{n}_{s}),
$$
which along with \eqref{5.8.5}  
implies that
\begin{equation}
                                           \label{1.17.5}
\kappa^{n\varepsilon}_{t\wedge\tau^{n}}:=
\hat{u}(x^{n}_{t\wedge\tau^{n}})e^{\phi^{n}_{t\wedge\tau^{n}}}
+\int_{0}^{t\wedge\tau^{n}}
f^{\alpha_{s}\beta^{n}_{s}   }
(p^{n}_{s},x^{n}_{s})
e^{-\phi^{n}_{s}}\,ds-\eta^{n\varepsilon}_{t}
 =\zeta^{n\varepsilon}_{t}+m^{n}_{t},
\end{equation}
where $\zeta^{n\varepsilon}_{t}$ is a decreasing process and
$$
\eta^{n\varepsilon}_{t}=\eta^{n\varepsilon}_{t}(\alpha_{\cdot})=
\varepsilon\delta_{1}^{-1}\int_{0}^{t\wedge\tau^{n}}e^{-\phi^{n}_{s}}
  \,ds+\int_{0}^{t\wedge\tau^{n}}e^{-\phi^{n}_{s}}[
|\xi^{n\varepsilon}
_{s}|+
 \chi_{\varepsilon}(x^{n}_{s}-x^{n}_{\kappa_{n}(s)})] \,ds.
$$

 Hence $\kappa^{n\varepsilon}_{t\wedge\tau^{n}}$ 
is at least a local supermartingale.  
 Assumption \ref{assumption 2.19.1}  
and  \eqref{2.29.2} show  that \eqref{2.29.1} and
the first inequality in \eqref{2.29.3}
hold.
It follows that the
 local supermartingale $\kappa^{n\varepsilon}_{t\wedge\tau^{n}}$  
 is, actually, a  supermartingale.

Furthermore,  obviously
$$
 \int_{0}^{ \infty}e^{-\phi^{n}_{s}}
 \,ds\leq 1/\chi  ,
$$
so that to prove the first assertion of the
 lemma, it only remains to  show that
\begin{equation}
                                           \label{1.17.6}
 \sup_{\alpha_{\cdot}\in\frA}E\int_{0}^{\infty}e^{-\phi^{n}_{s}}
\chi_{\varepsilon}(x^{n}_{s}-x^{n}_{\kappa_{n}(s)}) \,ds\to0 
\end{equation}
as $n\to\infty$. In light of the fact that
$c^{\alpha\beta}\ge\chi$, this is done in exactly the same way
as a similar fact is proved in \cite{Kr90}.  

That $\rho^{n\varepsilon}_{t\wedge\tau^{n}}(\alpha_{\cdot})$ is a local
supermartingale follows after computing its stochastic differential.
Then the fact that it is a supermartingale follows from 
the second estimate in \eqref{2.29.3} which is proved  by using
\begin{equation}
                                                      \label{6.12.2}
\int_{0}^{\infty}\lambda^{n}_{t}e^{-\psi_{t}^{n}}\,dt\leq1
\end{equation}
and the same argument as above. 
The lemma is proved.

{\bf Proof of Theorem \ref{theorem 4.17.1}}.
First we fix $x\in\bR^{d}$,
$\alpha_{\cdot}\in\frA$, and $\varepsilon>0$,
 take $\beta^{n}_{\cdot}(\alpha_{\cdot})$
form Lemma \ref{lemma 1.16.1} and prove that  the $\frB$-valued
functions defined on $\frA$ by
 $\bbeta^{n}(\alpha_{\cdot})=\beta^{n}_{\cdot}(\alpha_{\cdot})$
belong to $\bB$. To do that observe that if \eqref{4.5.4} holds and $T\leq 1/n$, then (a.s.)
$\beta^{n0}_{t}(\alpha^{1}_{\cdot})=\beta^{n0}_{t}(\alpha^{2}_{\cdot})$
for almost all $t\leq T$. By definition also (a.s.)
$$
p^{\alpha^{1}_{\cdot}\beta_{\cdot}^{n0}(\alpha^{1}_{\cdot})}_{s}
=p^{\alpha^{2}_{\cdot}\beta_{\cdot}^{n0}(\alpha^{2}_{\cdot})}_{s}
\quad\text{for almost all}\quad s\leq T.
$$
By uniqueness of solutions of \eqref{5.11.1} (see Assumption
\ref{assumption 5.23.2}),
 the processes $x_{t}^{n1}$ found from \eqref{4.5.1}
for $\alpha_{\cdot}=\alpha^{1}_{\cdot}$
and for $\alpha_{\cdot}=\alpha^{2}_{\cdot}$ coincide (a.s.)
for all $t\leq T$.

If \eqref{4.5.4} holds and $1/n<T\leq 2/n$, then by the above
solutions of \eqref{4.5.1} for $\alpha_{\cdot}=\alpha^{1}_{\cdot}$
and for $\alpha_{\cdot}=\alpha^{2}_{\cdot}$ coincide (a.s.)
for $t=1/n$ and then (a.s.)
$\beta^{n1}_{t}(\alpha^{1}_{\cdot} )=
\beta^{n1}_{t}(\alpha^{2}_{\cdot} )$ not only for all $t<1/n$
but also
for all $t\geq1/n$, 
which implies that (a.s.)
$$
p^{\alpha^{1}_{\cdot}\beta_{\cdot}^{n1}(\alpha^{1}_{\cdot})}_{s}
=p^{\alpha^{2}_{\cdot}\beta_{\cdot}^{n1}(\alpha^{2}_{\cdot})}_{s}
\quad\text{for almost all}\quad s\leq T
$$
and again the processes $x_{t}^{n }$ found from \eqref{4.5.1}
for $\alpha_{\cdot}=\alpha^{1}_{\cdot}$
and for $\alpha_{\cdot}=\alpha^{2}_{\cdot}$ coincide (a.s.)
for all $t\leq T$.

By induction we get that if \eqref{4.5.4} holds for a $T\in(0,\infty)$
and we define
$k $ as the integer such that $k/n<T\leq(k+1)/n$, then (a.s.)
\begin{equation}
                                               \label{4.5.5}
\beta^{n}_{t}(\alpha^{1}_{\cdot})=\beta^{nk}_{t}(\alpha^{1}_{\cdot} )=
\beta^{nk}_{t}(\alpha^{2}_{\cdot} )=
\beta^{n}_{t}(\alpha^{2}_{\cdot})\quad\text{for all}\quad
t< (k+1)/n ,
\end{equation}
$$
p^{\alpha^{1}_{\cdot}\beta_{\cdot}^{nk}(\alpha^{1}_{\cdot})}_{s}
=p^{\alpha^{2}_{\cdot}\beta_{\cdot}^{nk}(\alpha^{2}_{\cdot})}_{s}
\quad\text{for almost all}\quad s\leq T
$$
and  the processes $x_{t}^{n }$ found from \eqref{4.5.1}
for $\alpha_{\cdot}=\alpha^{1}_{\cdot}$
and for $\alpha_{\cdot}=\alpha^{2}_{\cdot}$ coincide (a.s.)
for all $t\leq T$.
This means that $\bbeta^{n} \in\bB$ indeed.

Furthermore, by the supermartingale property of 
$\rho^{n\varepsilon}_{t}(\alpha)$,
 for any
stopping times $\gamma^{\alpha_{\cdot}\beta_{\cdot} }\leq 
\tau^{\alpha_{\cdot}\beta_{\cdot} x}$ defined for  
each $\alpha_{\cdot}\in\frA$ and $\beta_{\cdot}\in\frB$ we have
$$
\hat{u}(x)\geq
E_{x}^{\alpha_{\cdot}\bbeta^{n}(\alpha_{\cdot})}
\hat{u}(x_{\gamma})e^{-\phi_{\gamma}-\psi_{\gamma}}
-E 
\eta^{n\varepsilon}_{\gamma}
(\alpha_{\cdot})e^{-\psi_{\gamma}}
$$
$$
+E_{x}^{\alpha_{\cdot}\bbeta^{n}(\alpha_{\cdot})}
 \int_{0}^{\gamma}
[f(p_{t},x_{t})+\lambda_{t}\hat{u}(x_{t})-
\lambda^{n}_{t}\eta^{n\varepsilon}_{t}
(\alpha_{\cdot})e^{\phi_{t}}]e^{-\phi_{t}-\psi_{t}}\,dt.
$$

Also observe that
$$
E_{x}^{\alpha_{\cdot}\bbeta^{n}(\alpha_{\cdot})}
[\int_{0}^{\gamma}\lambda_{t}\eta^{n\varepsilon}_{t}(\alpha_{\cdot})e^{-\psi_{t}}
\,dt+\eta^{n\varepsilon}_{\gamma}(\alpha_{\cdot})e^{-\psi_{\gamma}}]
$$
$$
\leq E 
\sup_{t\leq\gamma}\eta^{n\varepsilon}_{t}(\alpha_{\cdot})
\leq E 
 \eta^{n\varepsilon}_{\tau}(\alpha_{\cdot}).
$$

It follows that
$$
\hat{u}(x)\geq
E_{x}^{\alpha_{\cdot}\bbeta^{n}(\alpha_{\cdot})}
\big[
 \int_{0}^{\gamma}
[f(p_{t},x_{t})+\lambda_{t}\hat{u}(x_{t}) ]e^{-\phi_{t}-\psi_{t}}\,dt
$$
$$
+\hat{u}(x_{\gamma})e^{-\phi_{\gamma}-\psi_{\gamma}}\big]
-E 
 \eta^{n\varepsilon}_{\tau}(\alpha_{\cdot}),
$$
which owing to \eqref{3.3.1} yields
$$
\hat{u}(x)\geq\nliminf_{n\to\infty}\sup_{\alpha_{\cdot}\in\frA}
E_{x}^{\alpha_{\cdot}\bbeta^{n} (\alpha_{\cdot})}
\big[
 \int_{0}^{\gamma}
[f(p_{t},x_{t})+\lambda_{t}\hat{u}(x_{t}) ]e^{-\phi_{t}-\psi_{t}}\,dt
$$
$$
+\hat{u}(x_{\gamma})e^{-\phi_{\gamma}-\psi_{\gamma}}\big]
-\varepsilon/(\delta_{1}\chi)-Nd_{\varepsilon}(x).
$$
 
In light of the arbitrariness of $\varepsilon$ we arrive
at \eqref{2.11.5} and the theorem is proved.

For treating $\check{u}$ we use the following result.

\begin{lemma}  
                                                  \label{lemma 1.17.2}
Let $H[\check{u}]\geq0$ (everywhere)
  in $D$. Then for any $x\in\bR^{d}$, $\bbeta\in\bB$, and $\varepsilon>0$,
there exist a sequence $\alpha^{n}_{\cdot} 
\in\frA$, $n=1,2,...$, and a sequence of increasing continuous
$\{\cF_{t}\}$-adapted processes $\eta^{n\varepsilon}_{t}(\bbeta)$
 with $\eta^{n\varepsilon}_{0}(\bbeta)=0$
 such that 
the processes
$$
\kappa^{n\varepsilon}_{t} :=
\check{u}(x^{n}_{t})
e^{-\phi^{n}_{t}}
+\eta^{n\varepsilon}_{t} (\bbeta )
+\int_{0}^{t}
f^{n}_{s}
(p^{n}_{s},x^{n}_{s})
e^{-\phi^{n}_{s}}\,ds,
$$
where
\begin{equation}
                                               \label{5.6.4}
(x^{n}_{t},
\phi^{n}_{t})
=(x_{t},\phi_{t})
^{\alpha^{n}_{\cdot}\bbeta (\alpha^{n}_{\cdot})x},\quad
f^{n}_{t}(p,x)=
f^{\alpha^{n}_{t}\bbeta_{t}(\alpha^{n}_{\cdot}) }(p,x),
\quad p^{n}_{t}=
p^{\alpha^{n}_{\cdot}\bbeta (\alpha^{n}_{\cdot})}_{t},
\end{equation}
are submartingales on $[0,\tau^{\alpha^{n}_{\cdot}
\bbeta (\alpha^{n}_{\cdot})
  x}]$ and
\begin{equation}
                                                   \label{5.5.6}
\sup_{n}E 
\eta^{n\varepsilon}_{\infty}( \bbeta )<\infty,
\end{equation}
\begin{equation}
                                                   \label{5.5.7}
\nlimsup_{n\to\infty}
 E 
\eta^{n\varepsilon}_{\tau}(\bbeta)\leq
  \varepsilon/(\delta_{1}\chi) +Nd_{\varepsilon}(x),
\end{equation}
where $\delta_{1}$ is taken from Assumption \ref{assumption 5.23.1}
(ii).

 Furthermore,
if for any $n$ we are given a nonnegative,
  progressively measurable process 
$\lambda^{n}_{t}\geq0$ having finite integrals over
finite time intervals (for any $\omega$),
then the processes
$$
\rho^{n\varepsilon}_{t} :=
\check{u}(x^{n}_{t})
e^{-\phi^{n}_{t}
-\psi^{n}_{t}}
-\eta^{n\varepsilon}_{t} (\bbeta)e^{-\psi^{n}_{t}}
$$
$$
+\int_{0}^{t}
\big[f^{n}_{t}  
(p^{n}_{s},x^{n}_{s})
+\lambda^{n}_{s}\check{u}
(x^{n}_{s})
-\lambda^{n}_{s}\eta^{n\varepsilon}_{s} (\bbeta)
e^{ \phi^{n}_{s}}\big]
e^{-\phi^{n}_{s}
-\psi^{n}_{s}}\,ds 
$$
are submartingales on 
$[0,\tau^{\alpha^{n}_{\cdot}\bbeta (\alpha^{n}_{\cdot})x}]$, 
where we use  notation \eqref{5.6.4} and
$
\psi^{n}_{t} 
$ is taken from \eqref{4.28.1}.
 
  Finally,
$$
\sup_{n}E 
\sup_{t\geq0}|
\kappa^{n\varepsilon}_{t\wedge\tau }|<\infty,\quad
 \sup_{n}E 
\sup_{t\geq0}|
\rho^{n\varepsilon}_{t\wedge\tau }|<\infty.
$$
\end{lemma}  

Proof. Owing to Assumptions
\ref{assumption 1.9.1} and \ref{assumption 5.23.1} (i)
the function
$$
h(\alpha,x):=\inf_{\beta\in B}\big[ \bar L^{\alpha\beta}\check{u}(x)+
\bar f^{\alpha\beta}(x)\big]
$$
is a finite Borel function of $x$ and is continuous with 
respect to $\alpha$.
Its $\sup$ over $A$ can be replaced with the $\sup$
over an appropriate countable subset of $A$ and since
$$
\sup_{\alpha\in A}h(\alpha,x)\geq0,
$$
similarly to how $\beta(\alpha,x)$ was defined
in the proof of Lemma \ref{lemma 1.16.1}, one can find 
 a Borel function $\bar{\alpha}(x)$ in such a way that
\begin{equation}
                                                \label{1.17.07}
\inf_{\beta\in B}\big[\bar  L^{\bar{\alpha}(x)\beta}\check{u}(x)+
\bar f^{\bar{\alpha}(x)\beta}(x)\big]\geq  -\varepsilon 
\end{equation}
in $D$. If $x\not\in D$ we set $\bar{\alpha}(x)=\alpha^{*}$,
where $\alpha^{*}$ is a fixed element of $A$.  

After that we need some processes which we introduce recursively.
Fix $x$ and set $\alpha^{n0}_{t}\equiv\bar{\alpha}(x )$. Then define
$x^{n0}_{t}$ as a unique solution of the equation 
$$
x _{t}=x+\int_{0}^{t}\sigma(\alpha^{n0}_{s},
\bbeta _{s}(\alpha^{n0}_{\cdot }),p_{s}^{\alpha^{n0}_{\cdot}
\bbeta  (\alpha^{n0}_{\cdot })},x _{s})\,dw_{s} 
$$
$$
+\int_{0}^{t}b(\alpha^{n0}_{s},
\bbeta _{s}(\alpha^{n0}_{\cdot }),p_{s}^{\alpha^{n0}_{\cdot}
\bbeta  (\alpha^{n0}_{\cdot })},x _{s})\,ds.
$$
For $k\geq1$ introduce $\alpha^{nk}_{t}$ so that
$$
\alpha^{nk}_{t}=\alpha^{n(k-1)}_{t}\quad\text{for}\quad t<k/n,
$$
$$
\alpha^{nk}_{t}=\bar{\alpha}(x^{n(k-1)}_{k/n})
\quad\text{for}\quad t\geq k/n,
$$
where $x^{n(k-1)}_{t}$ is a unique solution of
 $$
x _{t}=x+\int_{0}^{t}\sigma(\alpha^{n(k-1)}_{s},
\bbeta _{s}(\alpha^{n(k-1)}_{\cdot }),p_{s}^{\alpha^{n(k-1)}_{\cdot}
\bbeta  (\alpha^{n(k-1)}_{\cdot })},x _{s})\,dw_{s} 
$$
\begin{equation}
                                                      \label{5.5.1}
+\int_{0}^{t}b(\alpha^{n(k-1)}_{s},
\bbeta _{s}(\alpha^{n(k-1)}_{\cdot }),p_{s}^{\alpha^{n(k-1)}_{\cdot}
\bbeta  (\alpha^{n(k-1)}_{\cdot })},x _{s})\,ds.
\end{equation}

As in the proof of Lemma \ref{lemma 1.16.1} we 
show that the above definitions
make sense as well as the definition
\begin{equation}
                                                \label{5.5.2}
\alpha^{n}_{t}=\alpha^{n(k-1)}_{t}\quad\text{for}\quad t<k/n.
\end{equation}
Next, by definition $x^{n}_{t}=
 x_{t}^{\alpha_{\cdot}^{n}\bbeta(\alpha^{n}_{\cdot})x}$ satisfies
$$
x_{t}=x+\int_{0}^{t}\sigma(\alpha^{n}_{s},\bbeta_{s}(
\alpha^{n}_{\cdot}),p_{s}^{\alpha^{n}_{\cdot} \bbeta (
\alpha^{n}_{\cdot})},x_{s})\,dw_{s}
$$
$$
+\int_{0}^{t}b(\alpha^{n}_{s},\bbeta_{s}(
\alpha^{n}_{\cdot}),p_{s}^{\alpha^{n}_{\cdot} \bbeta (
\alpha^{n}_{\cdot})},x_{s})\,ds.
$$
Equation  \eqref{5.5.2}   and the definitions
of $\bB$ and of control adapted processes show that
$x^{n}_{t}$ satisfies  \eqref{5.5.1} for $t\leq k/n$.
Hence, (a.s.) $x^{n}_{t}=x^{n(k-1)}_{t}$ for all $t\leq k/n$
and (a.s.) for all $t\geq0$, $\alpha^{n}_{t}
=\bar{\alpha}(x^{n}_{\kappa_{n}(t)})$ and 
$$
x^{n}_{t}=x+\int_{0}^{t}\sigma(\bar{\alpha}(x^{n}_{\kappa_{n}(s)}),
\bbeta _{s}(\alpha^{n}_{\cdot }),p^{n}_{s},x^{n}_{s})\,dw_{s}
$$
$$
+\int_{0}^{t}b(\bar{\alpha}(x^{n}_{\kappa_{n}(s)}),
\bbeta _{s}(\alpha^{n}_{\cdot }),p^{n}_{s},x^{n}_{s})\,ds,
$$
where $p^{n}_{s} =p_{s}^{\alpha^{n}_{\cdot} \bbeta (
\alpha^{n}_{\cdot})}$.

Now, 
introduce $\tau^{n}$ as the first exit time
of $x^{n}_{t}$ from $D$, set
$$
\beta^{n}_{s}=\bbeta_{s}(\alpha^{n}_{\cdot}),\quad
\phi^{n}_{t}=\phi^{\alpha^{n}_{\cdot}
\beta^{n}_{\cdot} x}_{t},\quad r^{n}_{s}=
r^{\alpha^{n}_{\cdot}
\beta^{n}_{\cdot} }(p^{n}_{s},x^{n}_{s}),
$$
where $r^{\alpha\beta}(p,x)$ is taken from 
Assumption  \ref{assumption 5.23.1} (ii),
 and observe that
by It\^o's formula  

$$
\check{u}(x^{n}_{t\wedge\tau^{n}})e^{-\phi^{n}_{t\wedge \tau^{n}}}
=\check{u}(x)+\int_{0}^{t\wedge\tau}e^{-\phi^{n}_{s}}
r^{n}_{s}\bar L^{\alpha^{n}_{s} 
\beta^{n}_{s}}
\check{u}(x^{n}_{s})
\,ds+m^{n}_{t},
$$
where $m^{n}_{s}$ is a martingale  and,
for $s<\tau^{n}$,
$$
\bar L^{\alpha^{n}_{s}
\beta^{n}_{s}}
\check{u}(x^{n}_{s})=\bar 
a_{ij}(\bar{\alpha}(x^{n}_{\kappa_{n}(s)}),\beta^{n}_{s}
,x^{n}_{s})
D_{ij}\check{u}(x^{n}_{s})
$$
$$
+\bar b_{i}(\bar{\alpha}(x^{n}_{\kappa_{n}(s)}),\beta^{n}_{s}
,x^{n}_{s})
D_{i}\check{u}(x^{n}_{s})
-\bar c(\bar{\alpha}(x^{n}_{\kappa_{n}(s)}),\beta^{n}_{s}
,x^{n}_{s})
 \check{u}(x^{n}_{s}).
$$
Similarly to the proof of Lemma \ref{lemma 1.16.1} we derive from 
\eqref{1.17.07} that,
for $s<\tau^{n}$,
$$
\bar L^{\alpha^{n}_{s}\beta^{n}_{s}}
\check{u}(x^{n}_{s})\geq -\varepsilon-\chi(x^{n}_{s}-x^{n}_{\kappa_{n}(s)}) -
\xi^{n\varepsilon}_{t}-
\bar 
f(\bar{\alpha}(x^{n}_{\kappa_{n}(s)}),\beta^{n}_{s},x^{n}_{\kappa_{n}(s)})
$$
$$
= -\varepsilon-\chi(x^{n}_{s}-x^{n}_{\kappa_{n}(s)}) -\xi^{n\varepsilon}_{t}
 -\bar f^{\alpha^{n}_{s}\beta^{n}_{s}}
( x^{n}_{s}),
$$
where $\xi^{n\varepsilon}_{t}$ are nonnegative
progressively measurable  processes satisfying \eqref{2.29.2}
and $\chi_{\varepsilon}(y)$ 
is a (nonrandom) bounded function on $\bR^{d}$ such that
$\chi_{\varepsilon}(y)\to0$ as $y\to0$.
It follows that
\begin{equation}
                                           \label{1.17.8}
\check{u}(x^{n}_{t\wedge\tau^{n}})e^{-\phi^{n}_{t\wedge\tau^{n}}}
+\int_{0}^{t\wedge\tau^{n}}
f^{\alpha^{n}_{s}\beta^{n}_{s}   }
(p^{n}_{s},x^{n}_{s})
e^{-\phi^{n}_{s}}\,ds+\eta^{n}_{t}
 =\zeta_{t}+m^{n}_{t},
\end{equation}
where $\zeta_{t}$ is an  increasing process and
$$
\eta^{n}_{t}=\eta^{n}_{t}(\bbeta)=
\varepsilon\delta_{1}^{-1}\int_{0}^{t\wedge\tau^{n}}e^{-\phi^{n}_{s}}
  \,ds+\int_{0}^{t\wedge\tau^{n}}e^{-\phi^{n}_{s}}[\xi^{n\varepsilon}
_{s}+
 \chi_{\varepsilon}(x^{n}_{s}-x^{n}_{\kappa_{n}(s)})] \,ds.
$$
 Hence the left-hand side of \eqref{1.17.8}
is a local submartingale and we finish the proof
in the same way as the proof of Lemma \ref{lemma 1.16.1}.
 The lemma is proved.

{\bf Proof of Theorem \ref{theorem 4.17.2}}.
 Similarly to the proof of Theorem \ref{theorem 4.17.1},
for any $\bbeta\in\bB$,
$$
\check{u}(x)\leq E_{x}^{\alpha^{n}_{\cdot}\bbeta (\alpha^{n}_{\cdot})}
\big[\int_{0}^{\gamma}[f(p_{t},x_{t})+\lambda_{t}\check{u}(x_{t})
+\lambda_{t}\eta_{t}^{n\varepsilon}(\bbeta)e^{ \phi_{t}}]e^{-\phi_{t}-\psi_{t}}\,dt
$$
$$
+\check{u}(x_{\gamma})e^{\phi_{\gamma}-\psi_{\gamma}}
+\eta^{n\varepsilon}_{\gamma}(\bbeta)e^{ \psi_{\gamma}}\big]
$$
$$
\leq E_{x}^{\alpha^{n}_{\cdot}\bbeta (\alpha^{n}_{\cdot})}
\big[\int_{0}^{\gamma}[f(p_{t},x_{t})+\lambda_{t}\check{u}(x_{t})
 ]e^{-\phi_{t}-\psi_{t}}\,dt
$$
$$
+\check{u}(x_{\gamma})e^{\phi_{\gamma}-\psi_{\gamma}}\big]
+E 
\eta^{n\varepsilon}_{\tau}
(\bbeta)  .  
$$
It follows that
$$
\check{u}(x)\leq \sup_{\alpha_{\cdot}\in\frA}
E_{x}^{\alpha _{\cdot}\bbeta(\alpha _{\cdot})}
\big[\int_{0}^{\gamma}[f(p_{t},x_{t})+\lambda_{t}\check{u}(x_{t})
 ]e^{-\phi_{t}-\psi_{t}}\,dt
$$
$$
+\check{u}(x_{\gamma})e^{\phi_{\gamma}-\psi_{\gamma}}\big]
+\nlimsup_{n\to\infty}
E 
\eta^{n\varepsilon}_{\tau}
(\bbeta)  ,
$$
$$
\check{u}(x)\leq \sup_{\alpha_{\cdot}\in\frA}
E_{x}^{\alpha _{\cdot}\bbeta(\alpha _{\cdot})}
\big[\int_{0}^{\gamma}[f(p_{t},x_{t})+\lambda_{t}\check{u}(x_{t})
 ]e^{-\phi_{t}-\psi_{t}}\,dt
$$
$$
+\check{u}(x_{\gamma})e^{\phi_{\gamma}-\psi_{\gamma}}\big]
 +\varepsilon/(\delta_{1}\chi)+Nd_{\varepsilon} ,
$$
which in light of the arbitrariness of $\varepsilon$ and $\bbeta\in\bB$
 finally yields \eqref{2.11.05} and the theorem is proved.

\mysection{Versions of Theorems \protect\ref{theorem 2.11.1},
\protect\ref{theorem 4.17.1}, and \protect\ref{theorem 4.17.2}
for ``uniformly nondegenerate'' case}
                                         \label{section 2.27.3}

In the first result of this section   $D$
is not assumed to be bounded.

Let  $\hat{u},\check{u} \in W^{2}_{1,loc}( D)\cap
C(\bar{D})$
  be given  functions for which there exist
sequences   $\hat{u}_{n},\check{u}_{n}\in C^{2}(\bar{D})$, $n\geq1$,
of $p$-insensitive in $D$ functions
which for each $n$ have 
uniformly continuous second-order derivatives (if $D$
is unbounded)
and such that $\hat{u}_{n},\check{u}_{n}$
 converge to $\hat{u}$ and $\check{u}$,
respectively,
uniformly in $\bar{D}$. For a sufficiently regular function
$u$ we denote by
$Du$ its  gradient  and $D^{2}u$ its Hessian.
In case of $\hat{u},\check{u}$ we take and fix
any Borel measurable versions of their gradients 
and Hessians.

\begin{theorem}
                                            \label{theorem 2.27.2}
Suppose that Assumptions
\ref{assumption 1.9.1}, \ref{assumption 5.23.2},
\ref{assumption 5.23.1} (i), (ii), 
Assumption
\ref{assumption 3.19.1} (ii), and Assumption  \ref{assumption 2.19.1}
(iii), (v)  are satisfied. Also suppose that a stronger
than Assumption  \ref{assumption 2.19.1}
(iv) is satisfied, namely, for any $x$
\begin{equation}
                                              \label{3.27.2}
\sup_{(\alpha_{\cdot},\beta_{\cdot} )\in\frA\times\frB }
E^{\alpha_{\cdot}\beta_{\cdot}}_{x}
\int_{0}^{\tau}
\sup_{\alpha\in A,\beta\in B}(
|\bar c^{\alpha\beta} -c^{\alpha\beta} _{\varepsilon}|
+|\bar f^{\alpha\beta} -f^{\alpha\beta} _{\varepsilon}|)(x_{t})
e^{-\phi_{t}}\,dt\to0.
\end{equation}
as $\varepsilon\downarrow0$. Finally, assume that
for any   $x\in D$
 \begin{equation}
                                              \label{3.26.2}
\sup_{(\alpha_{\cdot},\beta_{\cdot} )\in\frA\times\frB }
E^{\alpha_{\cdot}\beta_{\cdot}}_{x}
\int_{0}^{\tau}\big(|D^{2}\hat{u}-D^{2}\hat{u}_{n}|
+|D \hat{u}-D \hat{u}_{n}|)(x_{t})e^{-\phi_{t}}\,dt\to0
\end{equation}
as $n\to\infty$   and the same is true
if we replace $\hat{u}$ with $\check{u}$.

Then the following holds true:

(i) If $H[\hat{u}]\leq0$  
  in $D$ (a.e.)
and $\hat{u}\geq g$ on $\partial D$ (if $D\ne\bR^{d}$),
 then $v\leq\hat{u}$ in $\bar{D}$ and
\eqref{2.11.5} holds for any $\lambda^{\alpha_{\cdot}\beta_{\cdot}x}_{t}$
and $\gamma^{\alpha_{\cdot}\beta_{\cdot}x}$
  as in Theorem \ref{theorem 1.14.1}.

(ii) If $H[\check{u}]\geq0$  
  in $D$ (a.e.)
and $\check{u}\leq g$ on $\partial D$ (if $D\ne\bR^{d}$), 
then $v\geq\check{u}$ in $\bar{D}$ and
 \eqref{2.11.05} holds for any $\lambda^{\alpha_{\cdot}\beta_{\cdot}x}_{t}$
and $\gamma^{\alpha_{\cdot}\beta_{\cdot}x}$
  as in Theorem \ref{theorem 1.14.1}.

(iii) If $\hat{u}$ and $\check{u}$ are as in (i) and (ii) and
  $\hat{u}=\check{u}$, then
all assertions of Theorem \ref{theorem 1.14.1} hold true. Moreover,
  $v=\hat{u}$.

\end{theorem}

Before we proceed with the proof we note the following.
\begin{remark}
                                               \label{remark 4.16.1}
For $x\in\bR^{d}$ and
 $u=( u',u'')$, where $u'=(u'_{0},u'_{1},...,u'_{d})\in\bR^{d+1}$
 and $u''$
is in the set $\cS$ of $d\times d$ symmetric matrices,
introduce
\begin{equation}
                                                       \label{4.16.1}
H(u,x)=\supinf_{\alpha\in A\,\,\beta\in B}
\big(\bar a^{\alpha\beta}_{ij}(x)u''_{ij}+\sum_{i\geq1}
\bar b^{\alpha\beta}_{i}(x)u'_{i}-\bar c^{\alpha\beta}(x)u'_{0}
+\bar f^{\alpha\beta}(x)\big).
\end{equation}
Owing to Assumption \ref{assumption 1.9.1} (i) one can replace
$A$ and $B$ with their countable everywhere dense subsets.
Then we see that $H(u,x)$ is a Borel function of $x$.

Also note that  
$$
|H(u,x)|\leq N(\sum_{i,j=1}^{d}|u''_{ij} |+
\sum_{i=1}^{d}|u'_{i} |)
+
|u'_{0} |\sup_{\alpha,\beta}\bar c^{\alpha\beta}(x)+
\sup_{\alpha,\beta}|\bar f^{\alpha\beta}(x)|,
$$
$$
|H(u,x)-H(v,x)|\leq 
|u'_{0}-v'_{0}|\sup_{\alpha,\beta}\bar c^{\alpha\beta}(x)
$$
\begin{equation}
                                                 \label{3.27.5}
 +N\big(\sum_{i,j=1}^{d}|u''_{ij}-v''_{ij}|+
\sum_{i=1}^{d}|u'_{i}-v'_{i}|\big),
\end{equation}
where $N$ is independent of $u,v,x$. In light of
 \eqref{3.27.4}
 the right-hand sides are finite,
which, in particular, implies that
$H(u,x)$ is a Borel function of $(u,x)$.

If, in addition, $\bar c^{\alpha\beta}(x)$ and $\bar f^{\alpha\beta}(x)$
are bounded and continuous with respect to $x$
uniformly with respect to $(\alpha,\beta)$, then
the inequality
$$
|H(u,x)-H(u,y)|\leq N |x-y|(\sum_{i,j=1}^{d}|u''_{ij} |+
\sum_{i=1}^{d}|u'_{i} |)
$$
$$
+|u_{0}'|
\sup_{\alpha,\beta}|\bar c^{\alpha\beta}(x)-\bar c^{\alpha\beta}(y)|+
\sup_{\alpha,\beta}|\bar f^{\alpha\beta}(x)-\bar f^{\alpha\beta}(y)|
$$
shows that $H(u,x)$ is a continuous function of $x$,
which along with \eqref{3.27.5} guarantees that
$H(u,x)$ is a continuous function of $(u,x)$.

\end{remark}

{\bf Proof of Theorem \ref{theorem 2.27.2}}. (i) Introduce 
$\hat{h}_{n}=H[\hat{u}_{n}]$,
$$
c^{\alpha\beta}_{n}(p,x)=c^{\alpha\beta}(p,x)+n^{-1}
r^{\alpha\beta}(p,x),
$$
$$
f^{\alpha\beta}_{n}(p,x)= f^{\alpha\beta} (p,x)
 -r^{\alpha\beta}(p,x)\hat{h}_{n}(x)+n^{-1}
r^{\alpha\beta}(p,x)\hat{u}_{n}(x)
$$
$$
=r^{\alpha\beta}(p,x)[ \bar f^{\alpha\beta} ( x)
 -\hat{h}_{n}(x)+n^{-1}
 \hat{u}_{n}(x)],
$$ 

$$
L^{\alpha\beta}_{n}u(p,x)=L^{\alpha\beta}u(p,x)-
n^{-1}
r^{\alpha\beta}(p,x)u(x).
$$

Observe that $\hat{u}_{n}$ is $p$-insensitive
in $D$ with respect to $ L^{\alpha\beta}_{n}$. Owing to Definition
\ref{definition 5.8.1}, this follows from the fact that 
(dropping the superscripts $\alpha_{\cdot},\beta_{\cdot},x$)
for any $x\in D$ and $t<\tau$, we find that
 the coefficient of $dt$ in the stochastic
differential of
$$
\hat{u}_{n}(x_{t})e^{-\phi^{n}_{t}}
,\quad\text{where}\quad
\phi^{n}_{t}=\int_{0}^{t}c_{n}^{\alpha_{s}\beta_{s}}(p_{s},x_{s})\,ds,
$$
equals $e^{-\phi^{n}_{t}}$ times
$$
 L^{\alpha_{t}\beta_{t}}\hat{u}_{n}(p_{t},x_{t})-
n^{-1}
r^{\alpha_{t}\beta_{t}}(p_{t},x_{t})
\hat{u}_{n}(x_{t}) 
$$
$$
=r^{\alpha_{t}\beta_{t}}(p_{t},x_{t})
\big[\bar{L}^{\alpha_{t}\beta_{t}}\hat{u}_{n}(x_{t})
-n^{-1}\hat{u}_{n}(x_{t})
 \big]
=r^{\alpha_{t}\beta_{t}}(p_{t},x_{t}) 
L^{\alpha_{t}\beta_{t}}_{n}\hat{u}_{n}(\bar{p}_{t} ,x_{t}).
$$
Furthermore,
$$
\supinf_{\alpha\in A\,\,\beta\in B}
\big[
L^{\alpha\beta}_{n}\hat{u}_{n}(\bar{p},x)+
f^{\alpha\beta}_{n}(\bar{p},x)\big]=0,
$$
which makes us try to apply  Theorem \ref{theorem 4.17.1}
for each $n$.

Define $\hat{h}_{n\varepsilon}=H_{\varepsilon}[\hat{u}_{n}]$,
where $H_{\varepsilon} $ is constructed in the same way as
$H$ with $c_{\varepsilon}$ and $f_{\varepsilon}$
in place of $c$ and $f$, respectively, and
observe that, for each $n$ and $\varepsilon>0$,
 $\hat{h}_{n}$ is a Borel function on $\bar{D}$
and $\hat{h}_{n\varepsilon}$ is bounded and uniformly
continuous in $\bar{D}$ (cf.~Remark \ref{remark 4.16.1}). 
Also in $D$
$$
|\hat{h}_{n\varepsilon}-\hat{h}_{n }|=
|H_{\varepsilon}[\hat{u}_{n}]-H [\hat{u}_{n}]|
\leq (1+\sup|\hat{u}_{n}|) 
\sup_{\alpha\in A,\beta\in B}(
|\bar{c}^{\alpha\beta}-c^{\alpha\beta}_{\varepsilon}|
+|\bar{f}^{\alpha\beta}-f^{\alpha\beta}_{\varepsilon}|).
$$
Therefore, for any $x$  
$$
\lim_{\varepsilon\downarrow0}
\sup_{(\alpha_{\cdot},\beta_{\cdot} )\in\frA\times\frB }
E^{\alpha_{\cdot}\beta_{\cdot}}_{x}
\int_{0}^{\tau} |\hat{h}_{n\varepsilon}-\hat{h}_{n }|(x_{t})
e^{-\phi_{t}}\,dt=0.
$$

All other 
 assumptions of Theorem \ref{theorem 4.17.1}
are satisfied in light of the assumptions of the present theorem
and the fact that we added $n^{-1}$ to $\bar{c}$.
By Theorem \ref{theorem 4.17.1}  after setting
$$
\zeta^{\alpha_{\cdot}\beta_{\cdot}x}_{t}=
\int_{0}^{t}r^{\alpha_{s}\beta_{s}}(p^{\alpha_{\cdot}
\beta_{\cdot}}_{s},x^{\alpha_{\cdot}
\beta_{\cdot}x}_{s})\,ds
$$
we obtain
$$
\hat{u}_{n}(x)\geq\infsup_{\bbeta\in\bB\,\,\alpha_{\cdot}\in\frA}
E_{x}^{\alpha_{\cdot}\bbeta(\alpha_{\cdot})}\big[
\hat{u}_{n}(x_{\gamma})e^{-\phi_{\gamma}-\psi_{\gamma}-
\zeta_{\gamma}/n}
$$
$$
+\int_{0}^{\gamma}
\{f_{n}(p_{t},x_{t})+\lambda_{t}\hat{u}_{n}(x_{t})\}e^{-\phi_{t}-\psi_{t}
-\zeta_{t}/n}\,dt \big].
$$
 
Now we note that by considering $G+1$ in place of
$G$ we may assume that $G\geq1$ on $D$. We set $G:=1$ outside $D$.
Then, as follows easily from
 It\^o's formula,  the process 
$$
G(x_{t\wedge\tau})e^{-\phi_{t\wedge\tau}-\psi_{t\wedge\tau}}
+\int_{0}^{t\wedge\tau}( 1
+\lambda_{s}  ) 
e^{-\phi_{s}-\psi_{s}}\,ds
$$
is at least a local supermartingale, where
$$
(x_{t},\phi_{t},\tau)=(x_{t},\phi_{t},\tau)
^{\alpha_{\cdot}\beta_{\cdot}x},\quad
\psi_{t}=\psi_{t}^{\alpha_{\cdot}\beta_{\cdot}}
$$
and $x$, $\alpha_{\cdot}$, and $\beta_{\cdot}$ are arbitrary.
Nonnegative local supermartingales are supermartingales.
Therefore,
$$
E^{\alpha_{\cdot}\beta_{\cdot}}_{x}
e^{-\phi_{\gamma}-\psi_{\gamma}}+
E^{\alpha_{\cdot}\beta_{\cdot}}_{x}
\int_{0}^{ \gamma}( 1
+\lambda_{s}  ) 
e^{-\phi_{s}-\psi_{s}}\,ds\leq G(x).
$$

This shows that
$$
\hat{u}_{n}(x)
\geq
\infsup_{\bbeta\in\bB\,\,\alpha_{\cdot}\in\frA}
E_{x}^{\alpha_{\cdot}\bbeta(\alpha_{\cdot})}\big[
\hat{u} (x_{\gamma})e^{-\phi_{\gamma}-\psi_{\gamma}-\zeta_{t}/n}
$$
$$
+\int_{0}^{\gamma}
\{f (p_{t},x_{t})+\lambda_{t}\hat{u} (x_{t})\}e^{-\phi_{t}-\psi_{t}
-\zeta_{t}/n}\,dt \big]
$$
\begin{equation}
                                                 \label{3.26.5}
-(n^{-1}\delta_{1}^{-1}\sup_{D}|\hat{u}_{n}|+\sup_{D}|\hat{u}_{n}
-\hat{u}|) G(x)-\kappa_{n},
\end{equation}
where 
$$
\kappa_{n}(x)=\delta_{1}^{-1}\sup_{\alpha_{\cdot}\in\frA,
\beta_{\cdot}\in\frB}
E_{x}^{\alpha_{\cdot}\beta _{\cdot}}
\int_{0}^{\tau}(\hat{h}_{n}(x_{t}))^{+}e^{-\phi_{t}}\,dt.
$$
Observe that  
$$
(\hat{h}_{n})^{+}=(H[\hat{u}_{n}])^{+}
\leq(H[\hat{u}_{n}]-H[\hat{u}])^{+}
\leq N(|D^{2}(\hat{u}_{n}-\hat{u})|+
|D (\hat{u}_{n}-\hat{u})|)
$$
$$
+|\hat{u}_{n}-\hat{u}|
 \sup_{\alpha,\beta}\bar c^{\alpha\beta} .
$$
Here for any $\varepsilon>0$
$$
\sup_{\alpha,\beta}\bar c^{\alpha\beta}\leq
\sup_{\alpha,\beta}|\bar c^{\alpha\beta}-c^{\alpha\beta}_{\varepsilon}|
+\sup_{\alpha,\beta}|c^{\alpha\beta}_{\varepsilon}|,
$$
which along with our assumptions imply that
 $\kappa_{n}\to0$.
 After that by letting $n\to\infty$ in \eqref{3.26.5}
  we come to \eqref{2.11.5}.
Equation \eqref{2.11.5} with $\gamma=\tau$ and $\lambda\equiv0$
implies that $\hat{u}\geq v$.

(ii) Here the proof is very similar and yields \eqref{2.11.05},
from which we conclude that $\check{u}\leq v$.

(iii) The combination of assertions in (i) and (ii)
leads to $\hat{u}=\check{u}=v$, then \eqref{2.11.5}
and \eqref{2.11.05} imply that $v$ satisfies
\eqref{1.14.1} and,  since $\check{u}_{n}\to
\check{u}=v$ uniformly by assumption, $v$ is continuous
in $\bar{D}$ and in $\bR^{d}$. Finally, since $\hat{u}$ and 
$\check{u}$ have nothing to do with the fixed
control adapted process $p^{\alpha_{\cdot}\beta_{\cdot}}_{t}$,
the function $v$ is independent of the choice of this process.

The theorem is proved.

The assumptions of Theorem \ref{theorem 2.27.2}
admit an easy verification in the uniformly nondegenerate case.
 \begin{theorem}
                                            \label{theorem 2.27.3}
 Suppose that the domain $D$ is bounded,  
all requirements of  Assumptions 
\ref{assumption 1.9.1}, \ref{assumption 5.23.2}, and \ref{assumption 5.23.1} are satisfied, and 
$\hat{u}_{n}$ and $\check{u}_{n}$ not only converge
uniformly in $\bar{D}$ but also
converge in $W^{2}_{d}(D)$  
to $\hat{u}$ and $\check{u}$,
respectively. Then all assertions of Theorem
\ref{theorem 2.27.2} hold true.
\end{theorem}

Indeed, the existence of a
global barrier is well known for bounded domains
and uniformly nondegenerate operators with bounded
coefficients, so that Assumption
\ref{assumption 3.19.1} (ii) is satisfied.
Furthermore, in \eqref{3.27.2}
we can take $(c_{\varepsilon},f_{\varepsilon})
=(c^{(\varepsilon)},f^{(\varepsilon)})$
owing to Assumption \ref{assumption 5.23.1} (iii)  
and Lemma \ref{lemma 1.25.1}. The same lemma
guarantees that \eqref{3.26.2} holds and therefore
Theorem \ref{theorem 2.27.2} is applicable.

\begin{remark}
                                             \label{remark 4.29.1}
One may think that Theorem \ref{theorem 2.27.3}
is the only ``reasonable'' application of Theorem
\ref{theorem 2.27.2}. However, in a subsequent article we
will see an  application of Theorem
\ref{theorem 2.27.2} to a situation where $\check{u}$
depends only on part of the coordinates of a 
diffusion process, which does degenerate at each point,
but the above mentioned  part of it is uniformly nondegenerate.

\end{remark}

\mysection{An auxiliary result}
                                              \label{section 3.21.1}  

In this section $D$ is not assumed to be bounded.
We assume that we are given a continuous
$\cF_{t}$-adapted process  $x _{t}$  in $\bR^{d}$ and progressively measurable 
real-valued processes
$c _{t} $ and $f _{t} $. Suppose that $c_{t} \geq0$.

\begin{assumption}
                                        \label{assumption 3.21.1}
There exists a nonnegative bounded and
 continuous function $\Phi$ on $\bar{D}$
such that the process 
$$
\Phi(x  _{t\wedge\tau })e^{-\phi _{t\wedge\tau }}
+\int_{0}^{t\wedge\tau }|f_{s}|e^{-\phi _{s}}\,ds 
$$
is a supermartingale, where $\tau$ is the first exit time
of $x_{t}$ from $D$ and
$$
\phi_{t}=\int_{0}^{t}c_{s}\,ds.
$$

\end{assumption}

Let $D(n)$, $n\geq1$, be a sequence of
 subdomains of $D$. Introduce
$\tau_{n}$ as the first exit time
of $x_{t}$ from $D(n)$.

\begin{lemma}
                                          \label{lemma 3.19.1}
We have
\begin{equation}
                                               \label{6.12.1}
E\int_{0}^{\tau }
|f_{t}|e^{-\phi_{t}}\,dt 
\leq E\Phi(x_{0})I_{x_{0}\in D},
\end{equation}
\begin{equation}
                                               \label{3.21.2}
E\int_{\tau_{n} }^{\tau }
|f_{t}|e^{-\phi_{t}}\,dt 
\leq\sup_{(\partial D_{n})\setminus\partial D}\Phi
\quad(\sup_{\emptyset}:=0).
\end{equation}

\end{lemma}

Proof. By assumption, for any $t\in[0,\infty)$,
$$
E\big[
\Phi(x  _{t\wedge\tau })e^{-\phi _{t\wedge\tau }}
+\int_{0}^{t\wedge\tau }|f_{s}|e^{-\phi _{s}}\,ds \big]
\leq E\Phi(x_{0}),
$$
$$
E\int_{0}^{t\wedge\tau }|f_{s}|e^{-\phi _{s}}\,ds
\leq E[ \Phi(x_{0})-\Phi(x  _{t\wedge\tau })e^{-\phi _{t\wedge\tau }}]
\leq E\Phi(x_{0})I_{\tau>0},
$$
and sending $t\to\infty$ leads to \eqref{6.12.1}.

 Again by Assumption \ref{assumption 3.21.1}
for any $T\in[0,\infty)$
we have
$$
E\big[\Phi(x_{\tau_{n}\wedge T})e^{-\phi_{\tau_{n}\wedge T}}
+\int_{0}^{\tau_{n}\wedge T}|f_{t}|e^{-\phi_{t}}\,dt\big]
$$
$$
\geq
E\big[\Phi(x_{\tau \wedge T})e^{-\phi_{\tau \wedge T}}
+\int_{0}^{\tau \wedge T}|f_{t}|e^{-\phi_{t}}\,dt\big],
$$
$$
E\int_{\tau_{n}\wedge T}^{\tau\wedge T}
|f_{t}|e^{-\phi_{t}}\,dt\leq E
\big[\Phi(x_{\tau_{n}\wedge T})e^{-\phi_{\tau_{n}\wedge T}}
-\Phi(x_{\tau \wedge T})e^{-\phi_{\tau \wedge T}}\big]
$$
$$
=E
\big[\Phi(x_{\tau_{n}\wedge T})e^{-\phi_{\tau_{n}\wedge T}}
-\Phi(x_{\tau \wedge T})e^{-\phi_{\tau \wedge T}}\big]
I_{\tau_{n}<\tau  }
$$
$$
\leq E
 \Phi(x_{\tau_{n}\wedge T}) 
I_{\tau_{n}<\tau }.
$$
By sending $T\to\infty$ and using the monotone and dominated
convergence theorems we arrive at
$$
E\int_{\tau_{n} }^{\tau }
|f_{t}|e^{-\phi_{t}}\,dt\leq E
 \Phi(x_{\tau_{n} }) 
I_{\tau_{n}<\tau }
$$
and \eqref{3.21.2} follows. The lemma is proved.

\mysection{Proof of Theorem \protect\ref{theorem 1.14.1}}
                                       \label{section 4.25.1}

In this section all assumptions of Section \ref{section 2.26.3}
are supposed to be satisfied.

So far we did not use  Assumption \ref{assumption 3.19.1} (i)
concerning
the existence of $G$ vanishing on $\partial D$,
  which we need now to 
 deal with the case of general $D$. Take an expanding sequence  
of smooth domains $D_{n}\subset D$ from Assumption
\ref{assumption 3.22.1} 
and introduce
$$
v_{n}(x) 
=\infsup_{\bbeta\in\bB\,\,\alpha\in\frA}
E ^{\alpha_{\cdot}\bbeta(\alpha_{\cdot})}_{x}\big[
 \int_{0}^{\tau_{n}} f (p_{t},x_{t})e^{-\phi_{t}}\,dt
+g(x_{\tau_{n}})e^{-\phi_{\tau_{n}}}\big],
$$
where $\tau^{\alpha_{\cdot}\beta_{\cdot}x}_{n}$ is the first
exit time of $x^{\alpha_{\cdot}\beta_{\cdot}x}_{t}$
from $D(n)$. By Theorem \ref{theorem 2.27.3}
we have that $v_{n}$ are continuous in $\bR^{d}$ and
$$
v_{n}(x)=\infsup_{\bbeta\in\bB\,\,\alpha_{\cdot}\in\frA}
E_{x}^{\alpha_{\cdot}\bbeta(\alpha_{\cdot})}\big[
v_{n}(x_{\gamma_{n}})e^{-\phi_{\gamma_{n}}-\psi_{\gamma_{n}}}
$$
\begin{equation}
                                                     \label{7.10.1}
+\int_{0}^{\gamma_{n}}
\{f(p_{t},x_{t})+\lambda_{t}v_{n}(x_{t})\}e^{-\phi_{t}-\psi_{t}}\,dt \big],
\end{equation}
where
$ \gamma^{\alpha_{\cdot}\bbeta(\alpha_{\cdot})x}_{n}
=\gamma^{\alpha_{\cdot}\bbeta(\alpha_{\cdot})x}\wedge
\tau^{\alpha_{\cdot}\bbeta(\alpha_{\cdot})x}_{n}$.

We   claim that, as $n\to\infty$, 
\begin{equation}
                                               \label{3.20.7}
\kappa_{n}:=\sup_{\bR^{d}}|v_{n}-v|=
\sup_{D}|v_{n}-v|\to0 ,
\end{equation}
which, in particular, would imply the continuity of $v$
and the fact that $v$ is independent of the choice
of $p^{\alpha_{\cdot}\beta_{\cdot}}_{t}$.

To prove \eqref{3.20.7} introduce $v_{m}$ and $v_{nm}$
by replacing $g$ with $g_{m}$ in the definitions of $v$
and $v_{n}$, respectively, where $g_{m}$ are 
taken from the statement of the theorem. Observe that,
obviously, $\sup_{n}|v_{n}-v_{nm}|+|v-v_{m}|\to0$ as $m\to\infty$
uniformly on $\bR^{d}$. Therefore, {\em while proving\/}
\eqref{3.20.7} we may assume that  $\|g\|_{C^{2}(D)}<\infty$
and $g$ is $p$-insensitive.

 Then notice that in such a case we have
$$
E_{x}^{\alpha_{\cdot}\beta _{\cdot} }\big[\int_{0}^{ \tau}
 f (p_{t},x_{t})e^{-\phi_{t} }\,dt
+g(x_{ \tau})e^{-\phi_{\gamma\wedge\tau}
 }\big]
$$
$$
=g(x)+E_{x }^{\alpha _{\cdot} \beta 
 _{\cdot} } \int_{0}^{ \tau} 
 \hat{f} (p_{t},x_{t}) e^{-\phi_{t} }\,dt,
$$
where 
$$
\hat{f} ^{\alpha\beta}(p,x):= \big[f^{\alpha\beta} (p,x)
+r^{\alpha\beta}(p,x)\bar L^{\alpha\beta}g(x) \big]I_{D}(x) 
$$
satisfies  Assumption \ref{assumption 5.23.1} (i)-(iii).
Hence,
$$
u(x):=v(x)-g(x)=\infsup_{\bbeta\in\bB\,\,\alpha\in\frA}
E ^{\alpha_{\cdot}\bbeta(\alpha_{\cdot})}_{x}
 \int_{0}^{\tau }\hat f (p_{t},x_{t})e^{-\phi_{t}}\,dt.
$$

This argument   shows that we may also assume 
in the remaining part of the proof of \eqref{3.20.7} that $g=0$.
 Then
$$
v_{n}(x) 
=\infsup_{\bbeta\in\bB\,\,\alpha\in\frA}
E ^{\alpha_{\cdot}\bbeta(\alpha_{\cdot})}_{x}
 \int_{0}^{\tau_{n}} f (p_{t},x_{t})e^{-\phi_{t}}\,dt.
$$
 
Next, by 
using  It\^o's formula,
for any $\alpha_{\cdot}\in\frA$, $\beta_{\cdot}\in\frB$, and 
$x\in\bR^{d}$, we find that the process
\begin{equation}
                                                    \label{3.21.4}
G(x_{t\wedge\tau})e^{-\phi_{t\wedge\tau}}
+\int_{0}^{t\wedge\tau} 
e^{-\phi_{s}}\,ds
\end{equation}
is at least a local supermartingale, where
$$
(x_{t},\phi_{t},\tau)=(x_{t},\phi_{t},
\tau)^{\alpha_{\cdot}\beta_{\cdot}x}.
$$
Since $G$ is nonnegative in $D$,
the process \eqref{3.21.4} is a supermartingale (constant
if $x\not\in D$). 

Now, for $\chi>0$ introduce
$$
N_{\chi}=\sup_{(\alpha,\beta,x)\in A\times B\times D}
 |(\bar f^{\alpha\beta})^{(\chi)}|
$$
and observe that
by  Lemmas \ref{lemma 3.19.1} and \ref{lemma 1.25.1}
$$
|v_{n}(x)-v(x)|\leq I_{n}(x),
$$
where
$$
I_{n}(x):=\sup_{\alpha_{\cdot}\in\frA,\beta_{\cdot}
\in\frB }
E^{\alpha_{\cdot}\beta_{\cdot}}_{x}
\int_{\tau_{n}}^{\tau} |f (p_{t},x_{t})|e^{-\phi_{t}}\,dt
$$
$$
\leq \delta_{1}^{-1}\sup_{\alpha_{\cdot}\in\frA,\beta_{\cdot}
\in\frB }
E^{\alpha_{\cdot}\beta_{\cdot}}_{x}
\int_{\tau_{n}}^{\tau} |\bar f ( x_{t})|e^{-\phi_{t}}\,dt
$$
$$
\leq \delta_{1}^{-1}N_{\chi}\sup_{\alpha_{\cdot}\in\frA,\beta_{\cdot}
\in\frB }
E^{\alpha_{\cdot}\beta_{\cdot}}_{x}
\int_{\tau_{n}}^{\tau}  e^{-\phi_{t}}\,dt
$$
$$
+ \sup_{\alpha_{\cdot}\in\frA,\beta_{\cdot}
\in\frB }
E^{\alpha_{\cdot}\beta_{\cdot}}_{x}
\int_{0}^{\tau}|\bar f -\bar{f}^{(\chi)}|(x_{t}) e^{-\phi_{t}}\,dt
$$
$$
\leq \delta_{1}^{-1}N_{\chi}\sup_{(\partial D_{n})\setminus
\partial D}G
+N\|
\sup_{\alpha,\beta }
|\bar f^{\alpha,\beta} -(\bar f^{\alpha,\beta})^{(\chi)}|\,\|_{L_{d}(D)},
$$
where   $N$ is independent of $\chi,n$,
and $x$.
This and the fact that $G=0$ on $\partial D$
certainly imply \eqref{3.20.7}.

After that \eqref{7.10.1} (cf. \eqref{6.12.2})
yields
$$
v (x)\geq\infsup_{\bbeta\in\bB\,\,\alpha_{\cdot}\in\frA}
E_{x}^{\alpha_{\cdot}\bbeta(\alpha_{\cdot})}\big[
v (x_{\gamma_{n}})e^{-\phi_{\gamma_{n}}-\psi_{\gamma_{n}}}
$$
$$
+\int_{0}^{\gamma_{n}}
\{f(p_{t},x_{t})+\lambda_{t}v (x_{t})\}
e^{-\phi_{t}-\psi_{t}}\,dt \big]
 -\kappa_{n} .
$$

We use estimates like
$$
|v (x_{\gamma_{n}})e^{-\phi_{\gamma_{n}}-\psi_{\gamma_{n}}}
-v (x_{\gamma })e^{-\phi_{\gamma }-\psi_{\gamma }}|
=I_{\tau_{n}<\gamma}
|v (x_{\tau_{n}})e^{-\phi_{\tau_{n}}-\psi_{\tau_{n}}}
-v (x_{\gamma })e^{-\phi_{\gamma }-\psi_{\gamma }}|
$$
$$
\leq 2I_{\tau_{n}<\gamma}\sup_{t\in[\tau_{n},
\tau]}|v(x_{t})|e^{-\phi_{t}} ,
$$
$$
I_{\tau_{n}<\gamma}\int_{\gamma_{n}}^{\gamma}\lambda_{t}
|v(x_{t})|e^{-\psi_{t}}\,dt\leq
I_{\tau_{n}<\gamma}\sup_{t\in[\tau_{n},
\tau]}|v(x_{t})|e^{-\phi_{t}},
$$
where and sometimes in the future
  we drop the superscripts $\alpha_{\cdot}$,
$\beta_{\cdot}$, and $x$ for simplicity.

Then we see  that
$$
v (x)\geq\infsup_{\bbeta\in\bB\,\,\alpha_{\cdot}\in\frA}
E_{x}^{\alpha_{\cdot}\bbeta(\alpha_{\cdot})}\big[
v (x_{\gamma })e^{-\phi_{\gamma }-\psi_{\gamma }}
$$
\begin{equation}
                                                           \label{6.12.5}
+\int_{0}^{\gamma }
\{f(p_{t},x_{t})+\lambda_{t}v (x_{t})\}
e^{-\phi_{t}-\psi_{t}}\,dt \big]
 -\kappa_{n} -J_{n}(x),
\end{equation}
where 
$$
J_{n}(x)=I_{n}(x)+3R_{n}(x),
$$ 

$$
R_{n}(x):=\sup_{(\alpha_{\cdot},\beta_{\cdot})\in\frA\times\frB}
E^{\alpha_{\cdot} \beta_{\cdot}}_{x}I_{\tau_{n}<\gamma} 
 \sup_{t\in[\tau_{n},
\tau]}|v(x_{t})|e^{-\phi_{t}}. 
$$

To estimate $R_{n}(x)$  we observe that,
by Lemmas \ref{lemma 3.19.1} and \ref{lemma 1.25.1} for $x\in \bar{D}$
$$
|v(x)|\leq \delta_{1}^{-1}N_{\chi}G(x) +\varepsilon(\chi),
$$
where 
$$
\varepsilon(\chi)= N\|
\sup_{\alpha,\beta }
|\bar f^{\alpha,\beta} -(\bar f^{\alpha,\beta})^{(\chi)}|\,\|_{L_{d}(D)}
\to0,
$$
as $\chi\downarrow0$ (uniformly with respect to $x$).
Furthermore, since $G(x 
_{t\wedge\tau })\exp(-\phi_{t\wedge\tau })$  is a supermartingale,
 we have
$$
E^{\alpha_{\cdot} \beta_{\cdot}}_{x}  I_{\tau_{n}<\tau}
 \sup_{t\in[\tau_{n},
\tau]}[G(x_{t})e^{-\phi_{t}}]^{1/2}\leq N\big[
E^{\alpha_{\cdot} \beta_{\cdot}}_{x}  I_{\tau_{n}<\tau}
G(x_{\tau_{n}})e^{-\phi_{\tau_{n}}}\big]^{1/2},
$$
where $N$ is an absolute constant, and since  $v$ is bounded,

$$
E^{\alpha_{\cdot} \beta_{\cdot}}_{x} I_{\tau_{n}<\tau} 
 \sup_{t\in[\tau_{n},
\tau]}|v(x_{t})|e^{-\phi_{t}}\leq  N
E^{\alpha_{\cdot} \beta_{\cdot}}_{x}  I_{\tau_{n}<\tau}
 \sup_{t\in[\tau_{n},
\tau]}[|v(x_{t})|e^{-\phi_{t}}]^{1/2}
$$

$$
\leq NN_{\chi}\big[
E^{\alpha_{\cdot} \beta_{\cdot}}_{x}  I_{\tau_{n}<\tau}
G(x_{\tau_{n}}) \big]^{1/2}+N\varepsilon^{1/2}(\chi),
$$ 
where the constants $N$ are independent of $\chi$, $n$, and $x$.
By assumption $G=0$ on $\partial D$ and therefore we have that
$$
\sup_{(\alpha_{\cdot},\beta_{\cdot})\in\frA\times\frB}
E^{\alpha_{\cdot} \beta_{\cdot}}_{x} I_{\tau_{n}<\tau}
G(x_{\tau_{n}})\to0
$$
as $n\to\infty$ (uniformly with respect to $x$). It follows that
$$
\nlimsup_{n\to\infty}R_{n}(x)\leq N\varepsilon^{1/2}(\chi).
$$
Above we have also proved that
$$
\nlimsup_{n\to\infty}I_{n}(x)\leq N\varepsilon^{1/2}(\chi).
$$
It follows now from \eqref{6.12.5} that
$$
v (x)\geq\infsup_{\bbeta\in\bB\,\,\alpha_{\cdot}\in\frA}
E_{x}^{\alpha_{\cdot}\bbeta(\alpha_{\cdot})}\big[
v (x_{\gamma })e^{-\phi_{\gamma }-\psi_{\gamma }}
$$
$$
+\int_{0}^{\gamma }
\{f(p_{t},x_{t})+\lambda_{t}v (x_{t})\}
e^{-\phi_{t}-\psi_{t}}\,dt \big]
-  N\varepsilon^{1/2}(\chi),
$$
which after sending $\chi\downarrow0$ finally shows that
 $v(x)$ is greater than the right-hand side
of \eqref{1.14.1}. The reader understands that the opposite inequality
is proved similarly and this   brings the proof of
the theorem to an end.

\end{document}